\documentclass[11pt]{article}
\usepackage[T1]{fontenc}
\usepackage{lmodern}
\usepackage{amsmath,amssymb,amsthm,mathtools}
\usepackage{microtype}
\usepackage{cite}
\usepackage[hidelinks]{hyperref}
\usepackage[letterpaper,margin=1in]{geometry}
\hypersetup{
  pdftitle={Exact extremizers for judicious bisections of triangle-free graphs},
  pdfauthor={Shufei Wu},
  pdfkeywords={judicious partition, triangle-free graph, bisection, phase transition, extremal stability}
}
\allowdisplaybreaks[2]
\AtBeginDocument{%
\setlength{\abovedisplayskip}{6pt plus 2pt minus 3pt}
\setlength{\belowdisplayskip}{6pt plus 2pt minus 3pt}
\setlength{\abovedisplayshortskip}{3pt plus 2pt}
\setlength{\belowdisplayshortskip}{5pt plus 2pt minus 2pt}%
}

\newtheorem{theorem}{Theorem}[section]
\newtheorem{lemma}[theorem]{Lemma}
\newtheorem{proposition}[theorem]{Proposition}
\newtheorem{corollary}[theorem]{Corollary}
\newtheorem{remark}[theorem]{Remark}
\newtheorem{example}[theorem]{Example}

\newcommand{\Db}{D^{\mathrm b}}

\title{Exact extremizers for $\ell_p$-judicious bisections of triangle-free graphs}
\author{Shufei Wu\thanks{E-mail: \texttt{shufeiwu@hotmail.com}}\\
\small{School of Mathematics and Information Science, Henan Polytechnic University,
Henan, China, 454003}}
\date{}

\begin{document}
\maketitle

\begin{abstract}
Balogh, Clemen and Lidick\'y asked how large the
minimum $\ell_p$-norm of the internal-edge counts
over balanced partitions of an $n$-vertex
triangle-free graph can be. For each fixed
$2<p<\infty$ and all sufficiently large even $n$,
we determine this value exactly for bisections
and characterize all extremal graphs. They are
precisely the complete bipartite graphs
$K_{a,n-a}$ for which $a\le n/2$ maximizes an
explicit integer objective. All maximizing values
of $a$ lie in an interval of length bounded in
terms of $p$ alone.
We also prove quantitative stability for the
same range of $p$ and $n$: for sufficiently small
$\delta\ge 0$, every $n$-vertex triangle-free graph
with bisection value at least the asymptotic
optimum minus $\delta n^2$ is within
$O_p(\delta n^2+n)$ edge edits of a complete
bipartite graph with the optimal limiting
class ratio.
\end{abstract}

\medskip



\section{Introduction}
A partition of the vertex set of a graph into two classes
determines a cut, consisting of the edges joining the
classes. Maximizing the size of a cut is equivalent to
minimizing the total number of edges induced by the two
classes, but does not control how these edges are
distributed between them.
Judicious partition problems seek simultaneous control
of the internal-edge counts of the partition classes
\cite{BS1993,BS1999,BS2002,ABKS,BS2010,HouZeng2017,MaYu2016}.
This differs from the sparse-half problem
\cite{Erdos1997,Krivelevich,KeevashSudakov,NorinYepremyan,RazborovHalves},
which seeks a half inducing few edges without controlling
the other half. 
The $\ell_p$-norm provides a natural way
to combine the two internal-edge counts, interpolating
between their sum and their maximum
\cite{Scott,JMYU}.
Requiring the two classes to have equal size introduces
an additional difficulty, since any local adjustment
must preserve the balance of the partition.

Throughout, all graphs are finite, simple, and undirected. A
\emph{bisection} of a graph of even order is a partition $V(G)=A\cup B$
with $|A|=|B|$. For $X\subseteq V(G)$, let $e_G(X)$ denote the number
of edges induced by $X$, and for disjoint $X,Y\subseteq V(G)$, let
$e_G(X,Y)$ denote the number of edges between them. We omit the
subscript when the graph is clear. For $1\le p<\infty$, define
\[
\Db_p(G)=
\min_{\substack{V(G)=A\cup B\\|A|=|B|}}
\bigl(e(A)^p+e(B)^p\bigr)^{1/p}.
\]
For even $n$, let $M_p(n)$ be the maximum of $\Db_p(G)$ over all
triangle-free graphs on $n$ vertices. Determining $M_p(n)$ is the
bisection case of the norm problem for balanced partitions posed by
Balogh, Clemen and Lidick\'y \cite[Question~1.7]{BCL}.

The two endpoint objectives have different extremal complete bipartite graphs. Balogh, Clemen and Lidick\'y \cite{BCL} proved that, for sufficiently large even $n$, $M_1(n)\le n^2/16$, with equality only for $K_{n/4,3n/4}$. This bound extends to every $1\le p\le 2$, with $\Db_p(G)\le n^2/16$ for every triangle-free graph $G$ of sufficiently large even order $n$. Equality holds if and only if $4\mid n$ and $G\cong K_{n/4,3n/4}$ (see Proposition~\ref{rem:low-parameters}). The sum-norm objective can also be viewed as the minimum number of edges whose deletion makes the graph balanced bipartite. 
For $K_4$-free graphs, Balogh, Buczek, Grzesik and Kuc
\cite{BBGK2026} established the corresponding sharp bound
of $n^2/9$, resolving Conjecture~2.4 of Balogh, Clemen
and Lidick\'y \cite{BCL}.
At the maximum-norm endpoint, the sharp bound is $n^2/18$ when $n$ is sufficiently large and divisible by twelve, attained by $K_{n/3,2n/3}$ \cite{BCL}. Thus the optimal class ratio within the complete bipartite family depends on the norm. For fixed finite $p>2$, the remaining questions are to determine the optimal ratio and whether complete bipartite graphs account for all extremizers.

To state the main result, we introduce an integer optimization over
complete bipartite graphs. Let $N=n/2$. If a bisection of $K_{a,n-a}$,
where $0\le a\le N$, places $k$ vertices from the class of size $a$ in
one part, then the two parts induce $k(N-k)$ and $(a-k)(N-a+k)$ edges,
respectively. Define
\begin{equation}\label{eq:discrete-model}
\begin{split}
H_{p,n}(a)&=\min_{\substack{k\in\mathbb Z\\0\le k\le a}}
\left([k(N-k)]^p+[(a-k)(N-a+k)]^p\right)^{1/p},\\
B_p(n)&=\max_{\substack{a\in\mathbb Z\\0\le a\le N}}H_{p,n}(a),
\qquad
\mathcal A_p(n)=\{a\in\mathbb Z:0\le a\le N,\ H_{p,n}(a)=B_p(n)\}.
\end{split}
\end{equation}

\begin{theorem}\label{thm:exactification}
For every fixed $2<p<\infty$, there is $n_p$
such that, for every even $n\ge n_p$,
we have $M_p(n)=B_p(n)$.
Moreover, the extremal graphs, up to
isomorphism, are precisely
$K_{a,n-a}$ with $a\in \mathcal A_p(n)$.
\end{theorem}

The exact classification requires only a bounded number of candidate
class sizes. By Corollary~\ref{cor:integer-window}, all maximizing
values of $a$ lie in an interval of integers whose width depends only
on $p$; the integer optimization retains every possible tie.

The location of this interval is governed by the continuous relaxation
of the bipartite model, in which both the class ratio and the bisection
are allowed to vary. The resulting profile has a unique maximizing
smaller-class ratio $\sigma_p$ and a normalized maximum $\gamma_p$,
defined in Section~\ref{sec:model}. For $2<p<\infty$, we have
$1/4<\sigma_p<1/3$, with $\sigma_p\to1/4$ as $p\to2$ and
$\sigma_p\to1/3$ as $p\to\infty$. The corresponding values $\gamma_p$
decrease from $1/16$ to $1/18$. The profile decreases at least linearly
away from its maximum, and this loss, together with the rounding error,
yields the bounded integer window. In particular, as $n\to\infty$
through even integers,
\begin{equation}\label{eq:intro-asymptotic}
M_p(n)=\gamma_p n^2+O_p(n),\qquad
\Db_p(K_{\lfloor\sigma_p n\rfloor,n-\lfloor\sigma_p n\rfloor})
=\gamma_p n^2+O_p(n).
\end{equation}

The extremal classification is accompanied by a quantitative stability
result. A triangle-free graph with bisection value close to the
asymptotic optimum must be close to a complete bipartite graph with
nearly optimal class sizes. To state this precisely, define the edit
distance to $K_{a,n-a}$ by
\[
d_{\rm edit}(G,K_{a,n-a})=
\min_{\substack{V(G)=X\cup Y\\|X|=a}}
\bigl(e(X)+e(Y)+|X||Y|-e(X,Y)\bigr).
\]

\begin{theorem}\label{thm:intro-stability}
Fix $2<p<\infty$. There are positive constants
$C_p$ and $\delta_p$, and an integer $n_p$
such that, for every even $n\ge n_p$ and
every $0\le\delta\le\delta_p$, any triangle-free
graph $G$ on $n$ vertices satisfying
$D_p^{\mathrm b}(G)\ge(\gamma_p-\delta)n^2$
admits a partition $V(G)=X\cup Y$ with
\[
\frac{\bigl||X|-\sigma_p n\bigr|}{n}
+
\frac{e(X)+e(Y)+|X||Y|-e(X,Y)}{n^2}
\le C_p(\delta+n^{-1}).
\]
Moreover,
\[
d_{\mathrm{edit}}
\bigl(G,K_{\lfloor\sigma_p n\rfloor,
n-\lfloor\sigma_p n\rfloor}\bigr)
\le C_p(\delta+n^{-1})n^2.
\]
\end{theorem}

The continuous model also exhibits a transition in the structure of its minimizing bisections. At the optimal ratio, there are two competing types: one
leaves the smaller class intact, while the other splits it unevenly
for $2<p<p_c$ and evenly for $p\ge p_c$, where
$p_c\in(8/3,11/4)$ (see Theorem~\ref{thm:contact-phase}). This transition
concerns the minimizing bisections, rather than the complete bipartite
structure of the extremal graphs.

Related bisection problems have been studied under
various structural assumptions. Degree conditions
have played an important role in obtaining judicious
bisection bounds
\cite{XYY2010,XYYNote2010,XY2014}, while Lee, Loh
and Sudakov \cite{LLS2013} obtained further results
involving degree conditions and the structure of
tight components.
Under additional forbidden-subgraph and degree
restrictions, related estimates were obtained by
Fan, Hou and Yu \cite{FHY2018} and by Wu and Xiong
\cite{WX2024}. Lin and Zeng \cite{LinZeng2021}
studied maximum bisections of graphs without short
even cycles and, under additional cycle restrictions,
established a dichotomy between a large bisection
and a nearly bipartite structure.

The proof of the main results must overcome two distinct obstacles.
The first is to exclude graphs with $\alpha(G)<n/2$ from the
near-extremal regime, where $\alpha(G)$ denotes the independence
number. Estimates involving the maximum degree $\Delta(G)$ are useful
for this purpose, since every neighborhood in a triangle-free graph is
independent. Sampling bisections yields the required gap over much of
the parameter range, but near $p=2$ the leading-order estimates are
insufficient (see Example~\ref{ex:sampling-obstruction}). We overcome this
difficulty by exploiting the structure forced by the failure of these
estimates to construct a different bisection, while a density argument
handles the remaining cases. Once an independent half is available,
the bipartite profile constrains the class size, and sampling bounds
the remaining edge-count slack. Triangle-freeness then allows this
slack to be converted into an edit-distance bound.

The second obstacle is to pass from stability to exact classification.
The complete bipartite lower bound, together with the stability theorem,
shows that every extremal graph is within $O_p(n)$ edge edits of a
complete bipartite graph. This does not, however, rule out individual
missing or misplaced edges, and hence falls short of an exact
classification. Pikhurko, Slia\v{c}an and Tyros \cite{PST} developed
perfect-stability criteria for extremal problems whose objectives are
linear combinations of fixed-order subgraph densities. Since our
objective involves a minimum over bisections, we prove the required
zero-defect statement directly. A homogeneous estimate bounds the edge
defect in terms of a nonnegative integer defect, without an additive
error. Sampling bisections with the vertex counts prescribed by an
integer minimizer then shows that every positive defect incurs a strict
loss. The argument retains all ties in the integer optimization and
remains valid at the transition point $p=p_c$.

Section~\ref{sec:model} studies the bipartite profile and its contact
transition, while Section~\ref{sec:tools} develops the sampling
estimates. Sections~\ref{sec:unified}--\ref{sec:exactification} establish
the independence gap, stability and exact classification. The
appendices provide the supporting analytic estimates and exact
certificates, including the critical scaling analysis.

\section{The complete bipartite model}\label{sec:model}
We begin with the low-parameter range, where the
corresponding result follows from the theorem of
Balogh, Clemen and Lidick\'y \cite{BCL}.

\begin{proposition}\label{rem:low-parameters}
There exists $n_0$ such that, for every even
$n\ge n_0$, every $n$-vertex triangle-free
graph $G$, and every $1\le p\le 2$, we have
\[
D_p^{\mathrm b}(G)\le \frac{n^2}{16}.
\]
Equality holds if and only if $4\mid n$ and
$G\cong K_{n/4,3n/4}$.
\end{proposition}

\begin{proof}
Balogh, Clemen and Lidick\'y \cite[Theorem~1.1]{BCL} give
$\Db_1(G)\le n^2/16$, with equality only for $K_{n/4,3n/4}$.
Since $\|z\|_p\le\|z\|_1$ for $z\ge0$, minimizing over bisections gives
the upper bound, with equality possible only for the same graph.
Conversely, suppose $4\mid n$ and consider $K_{n/4,3n/4}$. If $xn$ vertices of its smaller class lie in one half, the
normalized class-edge counts are $a=x(1/2-x)$ and $b=1/16-x^2$, and
\[
 a^2+b^2-\frac1{256}=\frac{x^2(4x-1)^2}{8}\ge0.
\]
Thus $\|(a,b)\|_p\ge\|(a,b)\|_2\ge1/16$ for $p\le2$, while $x=0$ gives
equality. Hence this graph has $\Db_p=n^2/16$.
\end{proof}

\subsection{The complete bipartite profile}

Fix \(2<p<\infty\). For the complete bipartite graph with smaller-class
ratio \(s\), let \(\theta\) be the fraction of that class in the first
half and put \(t=1-2\theta\), after exchanging the halves if necessary.
The endpoints \(t=0,1\) represent an even split and an unsplit smaller
class. 
To analyze the minimum over bisections, we express
the $\ell_p$ unit circle in sum-and-difference coordinates.
Define \(U_p(t)\in[t,2]\), for \(0\le t\le1\), by
\begin{equation}\label{eq:Up}
 \left(\frac{U_p(t)-t}{2}\right)^p+
 \left(\frac{U_p(t)+t}{2}\right)^p=1.
\end{equation}
For fixed $t\in[0,1]$, the function
\[
 u\longmapsto
 \left(\frac{u-t}{2}\right)^p+
 \left(\frac{u+t}{2}\right)^p
\]
is strictly increasing on $[t,2]$, takes the value $t^p\le1$ at $u=t$,
and is strictly larger than $1$ at $u=2$. Thus $U_p(t)$ is well defined
and unique.

Set
\begin{equation}\label{eq:intro-parameters}
 R_p(t)=\frac{U_p(t)-1}{1-t^2},\qquad
 \kappa_p=\max_{0\le t\le1}R_p(t),\qquad
 \sigma_p=\frac{\kappa_p}{1+2\kappa_p},\qquad
 \gamma_p=\frac{\kappa_p}{2(1+2\kappa_p)^2}.
\end{equation}
Thus \(\kappa_p\) is the least parameter with
\(1+\kappa_p(1-t^2)\ge U_p(t)\) on \([0,1]\).
Equality defines a \emph{support contact}. The contact \(t=1\) is trivial;
for \(0\le t<1\), an additional contact is equivalently given by
\(R_p(t)=\kappa_p\).
We write \(\sigma_\infty=1/3\) and \(\gamma_\infty=1/18\) for the limiting
comparison constants as \(p\to\infty\); the central formula in
Corollary~\ref{cor:phase-explicit} gives these limits.
The convention \(R_p(1)=1/2\) is the continuous extension, since
\(U_p(t)-1\sim\frac{1-t^2}{2}\) as \(t\to1\).

\begin{lemma}
Extend $\kappa_p$, $\sigma_p$, and $\gamma_p$ in~\eqref{eq:intro-parameters}
to $p=2$ by
\[
 \kappa_2=\frac12,\qquad \sigma_2=\frac14,\qquad
 \gamma_2=\frac{\kappa_2}{2(1+2\kappa_2)^2}=\frac1{16}.
\]
Then $\kappa_p$, $\sigma_p$, and $\gamma_p$ are continuous in $p$ on
$[2,\infty)$. Moreover, for every $2<p<\infty$, one has
$\kappa_p>1/2$ and $\sigma_p>1/4$.
\end{lemma}

\begin{proof}
For $p=2$, equation~\eqref{eq:Up} gives
$U_2(t)=\sqrt{2-t^2}$. Hence
\[
 R_2(t)=\frac{\sqrt{2-t^2}-1}{1-t^2}
       =\frac{1}{\sqrt{2-t^2}+1},
\]
with the continuous value $R_2(1)=1/2$. Thus $\kappa_2=1/2$ and
$\sigma_2=1/4$.

For every fixed $P>2$, the function $R_p(t)$ extends jointly
continuously to $[2,P]\times[0,1]$. Indeed, implicit differentiation
of~\eqref{eq:Up} at $t=1$ gives $U_p(1)=1$ and $U_p'(1)=-1$.
Uniformly for $p$ in a compact interval $[2,P]$,
\(U_p(t)-1=\frac{1-t^2}{2}+o(1-t^2)\) as \(t\to1\).
Thus $R_p(t)\to1/2$ jointly as $t\to1$. Since the maximum is
taken over the compact interval $[0,1]$, the maximum theorem implies
that $\kappa_p$, and hence $\sigma_p$ and $\gamma_p$, are continuous
in $p$ on $[2,\infty)$.

For strictness, fix $2<p<\infty$, write $t=1-s$, and compare with
$u_0=1+(1-t^2)/2$.  As $s\to0$,
\[
 \left(\frac{u_0-t}{2}\right)^p+
 \left(\frac{u_0+t}{2}\right)^p
 =\left(s-\frac{s^2}{4}\right)^p+
 \left(1-\frac{s^2}{4}\right)^p
 =1-\frac p4s^2+o(s^2)<1,
\]
because $s^p=o(s^2)$.  The defining left side in~\eqref{eq:Up} is strictly
increasing in $u$, so $U_p(t)>u_0$ for all sufficiently small $s>0$ and
$R_p(t)>1/2$.
Therefore $\kappa_p>1/2$ and $\sigma_p>1/4$.
\end{proof}

For $0\le s\le1/2$ and $0\le\theta\le1$, define the coordinate maps and
the continuous complete bipartite profile by
\[
\begin{aligned}
 \mathcal C_1(s,\theta)&=\theta s\left(\frac12-\theta s\right),\\
 \mathcal C_2(s,\theta)&=(1-\theta)s
 \left(\frac12-(1-\theta)s\right),\\
 F_p(s)&=\min_{0\le\theta\le1}
 \bigl(\mathcal C_1(s,\theta)^p+\mathcal C_2(s,\theta)^p\bigr)^{1/p}.
\end{aligned}
\]

\begin{proposition}\label{prop:bipartite-profile}
The function $F_p$ is continuous on $[0,1/2]$ and satisfies
\(F_p(s)\le\gamma_p\), with equality if and only if $s=\sigma_p$.  On the entire right branch,
\begin{equation}\label{eq:profile-right-exact}
 F_p(s)=s(1/2-s)\qquad(\sigma_p\le s\le1/2).
\end{equation}
More quantitatively, if $t_p$ is any maximizing contact point of $R_p$ and
$\tau_p^-=1-t_p$, $\tau_p^+=1+t_p$, then
\begin{equation}\label{eq:linear-profile-corner}
 \gamma_p-F_p(s)\ge c_p^{\rm prof}|s-\sigma_p|
 \qquad(0\le s\le1/2),
\end{equation}
where
\[
\begin{aligned}
 c_p^{\rm prof}&=\min\{m_p^-,m_p^+\}>0,\qquad
 m_p^+=2\sigma_p-\frac12,\\
 m_p^-&=2^{-(p-1)/p}\frac{\tau_p^+}{4}
       (1-2\sigma_p\tau_p^+).
\end{aligned}
\]

For $a=\lfloor\sigma_pn\rfloor$ one has
\begin{equation}\label{eq:bipartite-rounding}
 \gamma_pn^2-2^{1/p-1}n
 \le \Db_p(K_{a,n-a})
 \le \gamma_pn^2+\frac n2,
\end{equation}
and
\[
 \max_{1\le b\le n/2}\Db_p(K_{b,n-b})
 =\gamma_pn^2+O_p(n).
\]
\end{proposition}

\begin{proof}
The coordinates $\mathcal C_1$ and $\mathcal C_2$ are precisely the normalized internal-edge
counts of a balanced partition of $K_{sn,(1-s)n}$ in which a proportion
$\theta$ of the smaller class is placed in the first half.  Thus $F_p(s)$ is
the continuous density-scale value of $n^{-2}\Db_p(K_{sn,(1-s)n})$.
Continuity follows from the minimum theorem on the compact rectangle
$[0,1/2]\times[0,1]$.

By symmetry we may take $0\le\theta\le1/2$.  For $0<s<1/2$, put
\[
 t=1-2\theta,\qquad E(s)=s\left(\frac12-s\right),
 \qquad \kappa_s=\frac{s}{1-2s}.
\]
Then the two coordinates, divided by $E(s)$, are
\[
 \frac{1+\kappa_s(1-t^2)-t}{2},\qquad
 \frac{1+\kappa_s(1-t^2)+t}{2}.
\]
For $s=\sigma_p$ we have $\kappa_s=\kappa_p$.  The definition of
$\kappa_p$ gives \(1+\kappa_p(1-t^2)\ge U_p(t)\)
for every $t\in[0,1]$, so every balanced partition at $s=\sigma_p$ has
$\ell_p$-cost at least $E(\sigma_p)=\gamma_p$.  The endpoint partition
$t=1$ attains this value.  At $s=1/2$, the two original vertex classes form
a balanced partition of cost zero, so the right-branch formula also holds
there.

If $\sigma_p\le s<1/2$, the endpoint partition gives $F_p(s)\le E(s)$.
Conversely $\kappa_s\ge\kappa_p$, so the normalized coordinates and the supporting
inequality show that every normalized vector has $\ell_p$-norm at least one.
Hence $F_p(s)=E(s)$, proving \eqref{eq:profile-right-exact}; since
$\sigma_p>1/4$, this is strictly below $E(\sigma_p)$ when $s>\sigma_p$.

Now let $0\le s<\sigma_p$, choose a maximizing contact $t_p$ of $R_p$, and put
\[
 a_s=\frac{s(1-t_p)}4\bigl(1-s(1-t_p)\bigr),\qquad
 b_s=\frac{s(1+t_p)}4\bigl(1-s(1+t_p)\bigr).
\]
Because $\kappa_p>R_p(1)=1/2$, we have $t_p<1$.  The contact equation gives
$U_p'(t_p)=-2\kappa_pt_p$ (also at $t_p=0$, by symmetry).  Writing
\[
 X=\frac{U_p(t_p)-t_p}{2},\qquad
 Y=\frac{U_p(t_p)+t_p}{2},
\]
implicit differentiation of $X^p+Y^p=1$ yields
\[
 U_p'(t_p)=-\frac{Y^{p-1}-X^{p-1}}{X^{p-1}+Y^{p-1}}>-1.
\]
Thus
\begin{equation}\label{eq:contact-slope}
 2\kappa_pt_p<1
 \quad\Longleftrightarrow\quad
 2\sigma_p(1+t_p)<1.
\end{equation}
Hence both $a_s$ and $b_s$ are strictly increasing on $[0,\sigma_p]$.
Put $g(s)=\|(a_s,b_s)\|_p$.  Since $b_s\ge a_s$ and
 $2\sigma_p(1+t_p)<1$, for $0<s\le\sigma_p$ we have
 \[
  a_s'\ge0,\qquad
  b_s'\ge\frac{\tau_p^+}{4}
       (1-2\sigma_p\tau_p^+),\qquad
  \frac{b_s}{g(s)}\ge2^{-1/p}.
 \]
 Therefore
 \[
  g'(s)\ge
  \left(\frac{b_s}{g(s)}\right)^{p-1}b_s'
  \ge m_p^-.
 \]
 The fixed-contact partition gives $F_p(s)\le g(s)$ and
 $g(\sigma_p)=\gamma_p$, so integration yields
 \[
  \gamma_p-F_p(s)\ge m_p^-(\sigma_p-s)
  \qquad(0\le s\le\sigma_p),
 \]
 with the endpoint $s=0$ obtained by continuity.  On the right branch,
 \eqref{eq:profile-right-exact} gives
 \[
  \gamma_p-F_p(s)
  =(s-\sigma_p)(s+\sigma_p-1/2)
  \ge m_p^+(s-\sigma_p).
 \]
 Since $\sigma_p>1/4$ and \eqref{eq:contact-slope} holds, both slopes are
 positive.  This proves \eqref{eq:linear-profile-corner}, and strict
 uniqueness follows at once.

For the discrete estimates, the two coordinates are uniformly
$1/2$-Lipschitz in $s$, hence $F_p$ is $2^{1/p-1}$-Lipschitz.  With
$a=\lfloor\sigma_pn\rfloor$ and $s=a/n$, every bisection of $K_{a,n-a}$
gives a feasible discrete value of $\theta$, so
\[
 \frac{\Db_p(K_{a,n-a})}{n^2}
 \ge F_p(s)\ge F_p(\sigma_p)-2^{1/p-1}|s-\sigma_p|.
\]
This proves the lower estimate in \eqref{eq:bipartite-rounding}; the endpoint
bisection and the $1/2$-Lipschitz bound for $E$ give the upper one.
Set $x=\theta s$. The coordinates $x(1/2-x)$ and
$(s-x)(1/2-s+x)$ are uniformly $1/2$-Lipschitz in $x$. Feasible values of
$x$ have spacing $1/n$, so rounding a minimizing continuous value changes
the normalized $\ell_p$-cost by $O_p(1/n)$, which gives the final $O_p(n)$ assertion.
\end{proof}

\begin{corollary}\label{cor:monotone}
The quantities $\kappa_p$ and $\sigma_p$ are strictly increasing in $p$
on $(2,\infty)$, while $\gamma_p$ is strictly decreasing in $p$ on
this interval.
\end{corollary}

\begin{proof}
Fix $0\le t<1$ and write
\[
 X_p=\frac{U_p(t)-t}{2},\qquad Y_p=\frac{U_p(t)+t}{2}.
\]
Then $X_p^p+Y_p^p=1$ and $0<X_p\le Y_p<1$ for $t<1$. Hence, if
$p_2>p_1>2$, then at $u=U_{p_1}(t)$,
\[
 \left(\frac{u-t}{2}\right)^{p_2}
 +\left(\frac{u+t}{2}\right)^{p_2}<1.
\]
The left side is strictly increasing in $u$, so
$U_{p_2}(t)>U_{p_1}(t)$ and therefore $R_{p_2}(t)>R_{p_1}(t)$ for
$0\le t<1$.
Let $t_1$ maximize $R_{p_1}$. Since
$\kappa_{p_1}>1/2=R_{p_1}(1)$, we have $t_1<1$, and thus
\[
 \kappa_{p_2}\ge R_{p_2}(t_1)>R_{p_1}(t_1)=\kappa_{p_1}.
\]
The formula $\sigma=\kappa/(1+2\kappa)$ is strictly increasing in
$\kappa$. Finally,
\[
 \frac{d}{d\kappa}\frac{\kappa}{2(1+2\kappa)^2}
 =\frac{1-2\kappa}{2(1+2\kappa)^3}<0
 \qquad(\kappa>1/2),
\]
which proves the assertion for $\gamma_p$.
\end{proof}
\subsection{The support contact transition}

We locate the maximizing contact of \(R_p\) using \(z=t^2\). The curvature
in this variable determines the monotonicity of the secant slopes.

\begin{lemma}\label{lem:curvature}
Let $W_p(z)=U_p(\sqrt z)$ on $[0,1]$.  If $2<p<3$, then $W_p$
has a unique inflection point $z_0\in(0,1)$, with
\[
 W_p''(z)<0\quad(0<z<z_0),\qquad
 W_p''(z)>0\quad(z_0<z<1).
\]
If $p\ge3$, then $W_p$ is convex on $[0,1]$.
\end{lemma}

\begin{proof}
Put $\alpha=p-1$.  Parametrize the rotated first-quadrant $\ell_p$ boundary by
\[
 u(\xi)=\frac{1+\xi}{(1+\xi^p)^{1/p}},\qquad
 z(\xi)=\frac{(1-\xi)^2}{(1+\xi^p)^{2/p}},
 \qquad 0\le\xi\le1.
\]
Then $W_p(z(\xi))=u(\xi)$ and
\begin{equation}\label{eq:Wpp-Q}
 W_p''(z(\xi))=
 \frac{(1+\xi^p)^{3/p}\xi^3Q_\alpha(\xi)}
 {4(1-\xi)^3(\xi+\xi^p)^3},
\end{equation}
where
\[
\begin{aligned}
 Q_\alpha(\xi)={}&1-\xi^{3\alpha}
 +\xi^{\alpha-1}\bigl((2\alpha+1)\xi-2\alpha\bigr)\\
 &+\xi^{2\alpha}\bigl(2\alpha\xi-(2\alpha+1)\bigr).
\end{aligned}
\]
The prefactor in \eqref{eq:Wpp-Q} is positive on $(0,1)$.

For $1<\alpha<2$, the function $Q_\alpha$ has exactly one zero
$\xi_Q\in(0,1)$ and changes sign there from positive to negative.  For
$\alpha\ge2$, one has $Q_\alpha\ge0$ on $[0,1]$.  The complete derivative
chain, including all endpoint values and the exact multiplicity at $\xi=1$,
is proved in Appendix~\ref{app:curvature-sign}.

For $1<\alpha<2$, logarithmic differentiation gives
\[
 \frac{z'(\xi)}{z(\xi)}
 =-\frac2{1-\xi}-\frac{2\xi^{p-1}}{1+\xi^p}<0.
\]
Thus $\xi:0\to1$ corresponds to $z:1\to0$.  Transferring the sign of
$Q_\alpha$ through \eqref{eq:Wpp-Q} and reversing the intervals proves the
first assertion, with $z_0=z(\xi_Q)$.  For $\alpha\ge2$, the second sign
statement and \eqref{eq:Wpp-Q} prove convexity.
\end{proof}

\begin{theorem}\label{thm:contact-phase}
There is a unique $p_c\in(2,3)$ satisfying
\begin{equation}\label{eq:pc-definition}
2^{2-1/p_c}\bigl(2^{1-1/p_c}-1\bigr)=p_c-1.
\end{equation}
If $2<p<p_c$, then $R_p$ has a unique maximizer
$t_p\in(0,1)$. If $p_c\le p<\infty$, then
$t=0$ is the unique maximizer.
\end{theorem}

\begin{proof}
For $0\le z<1$, define the secant slope and its numerator by
\[
 S_p(z)=\frac{W_p(z)-1}{1-z}=R_p(\sqrt z),\qquad
 \Phi_p(z)=(1-z)W_p'(z)+W_p(z)-1.
\]
Then
\begin{equation}\label{eq:SPhi-derivatives}
 S_p'(z)=\frac{\Phi_p(z)}{(1-z)^2},\qquad
 \Phi_p'(z)=(1-z)W_p''(z).
\end{equation}
Implicit differentiation of \eqref{eq:Up} at $t=1$ gives
$U_p(1)=1$ and $U_p'(1)=-1$. Since $z=t^2$, we have $W_p(1)=1$ and
$W_p'(1)=-1/2$. Consequently, $\Phi_p(1)=0$, and $S_p$ extends
continuously to $z=1$ with $S_p(1)=-W_p'(1)=1/2$.

Suppose first that $2<p<3$, and let $z_0$ be given by
Lemma~\ref{lem:curvature}.  Equation~\eqref{eq:SPhi-derivatives}
shows that $\Phi_p$ strictly decreases on $(0,z_0)$ and strictly increases
on $(z_0,1)$. Since the latter increase ends at $\Phi_p(1)=0$,
$\Phi_p(z)<0$ for $z_0\le z<1$.
If $\Phi_p(0)>0$, the first monotonicity interval contains exactly one zero
$z_p$, after which $\Phi_p$ stays negative.  Hence $S_p$ strictly increases
up to $z_p$ and strictly decreases afterwards, so $z_p$ is its unique global
maximizer.  If $\Phi_p(0)<0$, then $\Phi_p<0$ on $(0,1)$ and $S_p$ is
strictly decreasing.  The equality case is equally strict: if
$\Phi_p(0)=0$, then the first strict decrease makes $\Phi_p(z)<0$ for
$0<z\le z_0$, and the final strict increase reaches zero only at $z=1$.
Thus $z=0$ is uniquely maximizing whenever $\Phi_p(0)\le0$.

The sign at the center determines the remaining cases. Put $a_p=2^{1-1/p}$.
At $t=0$, implicit differentiation gives
\[
 W_p(0)=a_p,\qquad W_p'(0)=-\frac{p-1}{2a_p}.
\]
Therefore
\[
 \Phi_p(0)=a_p-1-\frac{p-1}{2a_p}
          =\frac{f(p)}{2a_p},
 \qquad
 f(p)=2a_p(a_p-1)-(p-1).
\]
Appendix~\ref{app:center-sign} verifies that $f$ is strictly decreasing on
$(2,\infty)$ and that $f(2)>0>f(3)$.  Thus $f$ has exactly one zero in
$(2,3)$, namely the number in
\eqref{eq:pc-definition}. Combining its sign with the preceding
classification proves the theorem for $2<p<3$, including the critical
equality case.

For $3\le p<\infty$, Lemma~\ref{lem:curvature} gives convexity of $W_p$ and hence
the endpoint-chord inequality
\[
 W_p(z)\le(1-z)W_p(0)+zW_p(1).
\]
It is strict for $0<z<1$: otherwise equality in convexity at an interior
point would make $W_p$ affine, whereas
\[
 W_p'(0)=-\frac{p-1}{2a_p}\ne-\frac12=W_p'(1)
\]
because $p-1\ge2>a_p$. Thus $S_p(z)<S_p(0)$ in the interior.
At the remaining endpoint, $S_p(1)=1/2<S_p(0)=a_p-1$ for
$3\le p<\infty$, since $a_p\ge\sqrt[3]{4}>3/2$. Thus the central contact is uniquely maximizing
also at and above the convexity threshold for $W_p$.

\end{proof}

Henceforth, for $2<p<\infty$, let $t_p$ denote the unique maximizer of
$R_p$; thus $t_p\in(0,1)$ for $2<p<p_c$ and $t_p=0$ for $p\ge p_c$.

\begin{corollary}
\label{cor:phase-explicit}
For every $p_c\le p<\infty$, $\kappa_p=2^{1-1/p}-1$, and hence
\[
 \sigma_p=\frac{2-2^{1/p}}{4-2^{1/p}},\qquad
 \gamma_p=
 \frac{2^{1/p}(2-2^{1/p})}{2(4-2^{1/p})^2}.
\]
\end{corollary}

\begin{proof}
Theorem~\ref{thm:contact-phase} gives the unique maximizing contact $t=0$,
so $\kappa_p=R_p(0)=U_p(0)-1=2^{1-1/p}-1$. Substitute this value into
\eqref{eq:intro-parameters}.
\end{proof}

\begin{remark}
The threshold $p=3$ concerns only the convexity of $W_p$: by then the
optimal contact is already central.  Thus $p=3$ creates neither a new
extremal phase nor a new graph-theoretic regime.  The central-contact
density argument applies throughout $p\ge p_c$.
\end{remark}

\section{Sampling tools}\label{sec:tools}
Throughout Sections~\ref{sec:tools}--\ref{sec:exactification},
$G$ denotes a triangle-free graph of even order $n$,
unless stated otherwise. The two sampling lemmas below
do not require triangle-freeness, with even order needed
only for their bisection conclusions.

\begin{lemma}\label{lem:balanced-sampling}
Let $U$ have size $u=\Theta(n)$, let $H_1,\dots,H_s$ be a fixed finite
family of graphs on $U$, and prescribe $|X|=\theta u+O(1)$.  There is a
choice of $X$ for which, simultaneously for all $i$,
\[
 e(H_i[X])=\theta^2e(H_i)+O(n^{3/2}).
\]
More generally, in a fixed-size random partition into
a fixed number of blocks, any fixed finite collection
of induced and mixed edge counts can be made
simultaneously within $O(n^{3/2})$ of their respective
expectations. The same conclusion holds when some sets
are fixed or when disjoint vertex sets are sampled
independently.
In particular,
if an $n$-vertex graph $G$ has $m$ edges, then it has a bisection
$V(G)=A\cup B$ with
\[
 e(A)=\frac m4+O(n^{3/2}),\qquad
 e(B)=\frac m4+O(n^{3/2}),
\]
and consequently, for every $p\in[1,\infty)$,
\[
 \Db_p(G)\le 2^{1/p}\frac m4+O(n^{3/2}).
\]
\end{lemma}

\begin{proof}
Write $|X|=k=\theta u+O(1)$.  For each edge of $H_i$, the probability that
both endpoints lie in a uniformly chosen $k$-set is $(k)_2/(u)_2$.  Hence
\[
 \mathbb E e(H_i[X])=\frac{(k)_2}{(u)_2}e(H_i)
 =\theta^2e(H_i)+O(n).
\]
The last error is uniform in $0\le\theta\le1$: the sampling ratio differs
from $\theta^2$ by $O(u^{-1})$, while $e(H_i)=O(u^2)$ and $u=\Theta(n)$.

Represent the prescribed sets by consecutive blocks of a uniform random
permutation, using independent permutations for disjoint vertex sets.
Expose these permutations one entry at a time. Changing one exposed
choice can be coupled by a single transposition, which changes each
relevant edge count by at most $Cn$. Its Doob martingale therefore has
$O(n)$ differences of size $O(n)$, so orthogonality gives
\(\operatorname{Var}Z=O(n^3)\) for every induced or mixed count $Z$ under consideration. Chebyshev's
inequality gives
\[
 \Pr\bigl(|Z-\mathbb EZ|>Kn^{3/2}\bigr)\le C/K^2.
\]
For sufficiently large fixed $K$, a union bound over the fixed finite
family leaves positive probability that all counts satisfy the required
estimate.

This argument includes complementary counts, edges to a fixed set, and
several prescribed blocks. Each expectation has the same $O(n)$
finite-size error; changing a prescribed size by $O(1)$ also changes a
count by at most $O(n)$. Thus the normalized errors are $O(n^{-1/2})$,
uniformly over the graph and all admissible proportions. For the final
assertion use two complementary blocks of size $n/2$ and include both
induced counts in the family.
\end{proof}

For finite $p\ge2$, coordinatewise concentration gives an $O(n^{3/2})$
error. At a fixed positive density-scale target, the following squared-norm
estimate yields an $O_p(n)$ error instead.

\begin{lemma}[Squared-norm sampling]\label{lem:squared-norm-sampling}
Fix \(2\le p<\infty\). For any nonnegative random vector \(Z\in\mathbb R^2\) with
finite second moment and any deterministic \(x\ge0\),
\begin{equation}\label{eq:reference-sampling}
 \mathbb E\|Z\|_p^2\le\|x\|_p^2+
 \nabla(\|\cdot\|_p^2)(x)\cdot(\mathbb EZ-x)
 +(p-1)\mathbb E\|Z-x\|_2^2.
\end{equation}
If \(Z\) consists of two edge counts in a fixed-size experiment as in
Lemma~\ref{lem:balanced-sampling}, allowing independent samples from
disjoint vertex sets, and \(\mu=\mathbb EZ\), then
\[
 \mathbb E\|Z\|_p^2\le\|\mu\|_p^2+
 (p-1)\mathbb E\|Z-\mu\|_2^2\le\|\mu\|_p^2+C_p n^3.
\]
In particular, if \(Z\) is the internal-edge vector of a random exact
bisection and \(\mathbb EZ=n^2z+O(n)\) coordinatewise, where \(z\ge0\)
is bounded, then
\begin{equation}\label{eq:squared-sampling}
 \left(\frac{\Db_p(G)}{n^2}\right)^2\le\|z\|_p^2+\frac{C_p}{n}.
\end{equation}
The sampling constants are uniform in the graph and proportions; they
may depend on the fixed number of blocks, and in the last assertion on
the bounds for \(z\) and the mean error.
\end{lemma}

\begin{proof}
Put $f_p(x)=\|x\|_p^2$ on the nonnegative quadrant. For $p>2$ and
$x_1,x_2>0$, set $S=x_1^p+x_2^p$ and
$v=(x_1^{p-1},x_2^{p-1})^{\mathsf T}$. Differentiation gives
\[
 \nabla^2 f_p(x)
 =2(p-1)S^{2/p-1}\operatorname{diag}(x_1^{p-2},x_2^{p-2})
 -2(p-2)S^{2/p-2}vv^{\mathsf T}
 \preceq 2(p-1)I_2.
\]
Indeed, $S^{2/p-1}x_i^{p-2}\le1$, and the subtracted matrix is positive
semidefinite. For positive coordinates the Hessian bound therefore gives
\[
 f_p(y)\le f_p(x)+\nabla f_p(x)\cdot(y-x)
 +(p-1)\|y-x\|_2^2.
\]
The general nonnegative case follows by approximation from the interior;
at the origin the squared norm is differentiable with gradient zero.
For $p=2$ the displayed inequality is an identity. Taking $y=Z$ and
expectations proves \eqref{eq:reference-sampling}. Taking $x=\mu$
cancels the linear term and proves the centered inequality.
For the fixed-size sampling assertion, the variance bound in the proof of
Lemma~\ref{lem:balanced-sampling} gives
$\mathbb E\|Z-\mu\|_2^2=\operatorname{Var}Z_1+\operatorname{Var}Z_2=O(n^3)$.
In the bisection
case, $\Db_p(G)^2\le\mathbb E\|Z\|_p^2$. Moreover,
$\|\mu\|_p^2=n^4\|z\|_p^2+O_p(n^3)$, which proves
\eqref{eq:squared-sampling}.
\end{proof}

\begin{lemma}\label{lem:alphaedges}
Every triangle-free graph $H$ satisfies
\[
 e(H)\le \alpha(H)(|V(H)|-\alpha(H)).
\]
\end{lemma}

\begin{proof}
Let $S$ be a maximum independent set and $T=V(H)\setminus S$. Since every
neighborhood is independent, $d(v)\le\alpha(H)$ for $v\in T$. Hence
\[
 e(S,T)+2e(T)=\sum_{v\in T}d(v)\le \alpha(H)|T|,
\]
and the result follows.
\end{proof}

\begin{lemma}\label{lem:fixedI}
Let $I$ be a maximum independent set of size $rn<n/2$ and put
$J=V(G)\setminus I$. Set
\(z=\frac{e(J)}{n^2}\) and \(w=\frac{e(I,J)}{n^2}\).
There is a bisection whose internal-edge densities are at most
\[
 \lambda^2z+\lambda w+O(n^{-1/2}),\qquad
 \mu^2z+O(n^{-1/2}),
\]
where \(\lambda=\frac{1/2-r}{1-r}\) and \(\mu=\frac1{2(1-r)}\).
If $e(G)>\rho n^2$, then the two main terms are at most
\begin{equation}\label{eq:ABr}
 A(r)=\frac{r(1/2-r)}{2(1-r)},\qquad
 B(r)=\frac{r(1-r)-\rho}{4(1-r)^2}.
\end{equation}
\end{lemma}

\begin{proof}
Choose uniformly a set $X\subseteq J$ of size $n/2-|I|$ and take
$I\cup X$ as one class. The proportions of $J$ in $X$ and
$J\setminus X$ are exactly $\lambda$ and $\mu$, respectively. Hence
\[
 \mathbb E\,e(I,X)=\lambda e(I,J)+O(n),
\]
\[
 \mathbb E\,e(X)=\lambda^2e(J)+O(n),\qquad
 \mathbb E\,e(J\setminus X)=\mu^2e(J)+O(n).
\]
By Lemma~\ref{lem:balanced-sampling}, applied to $X$, $J\setminus X$, and
the mixed edges between $I$ and $X$, there is a single choice of $X$ for
which all three variables differ from their expectations by $O(n^{3/2})$.
This proves the first assertion.

Since $\Delta(G)\le\alpha(G)=rn$,
\[
 z+w=\frac{e(G)}{n^2}\le\frac r2,
 \qquad
 w+2z=\frac1{n^2}\sum_{v\in J}d(v)\le r(1-r).
\]
Using $\lambda+\mu=1$, we obtain
\[
 \lambda^2z+\lambda w
 =\lambda\bigl((z+w)-\mu z\bigr)
 \le\frac{\lambda r}{2}
 =\frac{r(1/2-r)}{2(1-r)}.
\]
If $z+w>\rho$, then the second displayed inequality gives
$z<r(1-r)-\rho$, and consequently
\[
 \mu^2z<\frac{r(1-r)-\rho}{4(1-r)^2}.
\]
This gives \eqref{eq:ABr}, and hence internal-edge densities at most
$A(r)+O(n^{-1/2})$ and $B(r)+O(n^{-1/2})$, respectively.
\end{proof}

\section{Reduction to an independent half}\label{sec:unified}

Our aim is to exclude graphs without an independent half from the
near-extremal regime. We divide the argument according to whether the
optimal support contact is interior or central. In the interior range,
either of two canonical sampling estimates gives the required gap unless
their leading terms are both too large; in that case the resulting
structural constraints yield a different bisection by orienting defects
and restoring balance. In the central range, a density split selects
between a random bisection and a sample extending a maximum independent
set. Together these arguments give the uniform independence gap at the
end of the section.

\subsection{Interior contact: sampling or structure}

Set
\[
 c_*:=\frac1{16}-\frac1{5000}=\frac{623}{10000}.
\]
If $\alpha(G)<n/2$, the dichotomy below gives
$\Db_2(G)\le c_*n^2+O(n^{3/2})$. For $p>2$, norm monotonicity
implies $\Db_p(G)\le\Db_2(G)$. For $2<p<p_c$, $\gamma_p>\gamma_{p_c}>c_*$, where the latter
inequality is proved in Lemma~\ref{lem:phase-rational-margins}. Thus a fixed
positive gap absorbs the $O(n^{3/2})$ term for all sufficiently large $n$.

Throughout this subsection, let $I$ be a maximum independent set with
$|I|=rn<n/2$, and put $J=V(G)\setminus I$. Consider a random balanced
sample and a sample that keeps $I$ on one side and fills that side from
$J$. If both main-term costs exceed $c_*$, the profile constraints below
force the defect pattern used in the structural alternative.

Write
\[
 z=\frac{e(J)}{n^2},\qquad w=\frac{e(I,J)}{n^2},\qquad
 \lambda=\frac{1/2-r}{1-r},\qquad \mu=\frac1{2(1-r)}.
\]
\paragraph{Two canonical samples.}
The two main-term normalized $\ell_2$ costs are
\begin{equation}\label{eq:canonical-costs}
 \begin{aligned}
 C_{\rm rnd}&=\frac{z+w}{2\sqrt2},\\
 C_{\rm fix}&=\sqrt{(\lambda^2z+\lambda w)^2+(\mu^2z)^2}.
 \end{aligned}
\end{equation}
Indeed, a uniformly random bisection has main-term edge densities
$(z+w)/4$ on each side, so Lemma~\ref{lem:balanced-sampling} gives a
bisection of cost at most $C_{\rm rnd}n^2+O(n^{3/2})$.
For the fixed-$I$ sample, put all of $I$ in one side and choose uniformly
$(1/2-r)n$ vertices from $J$ to complete that side.
Lemma~\ref{lem:fixedI} gives the main coordinates
$\lambda^2z+\lambda w$ and $\mu^2z$, and hence cost at most
$C_{\rm fix}n^2+O(n^{3/2})$.  Both procedures preserve exact balance,
and the errors are uniform over the graphs under consideration.

\begin{lemma}\label{lem:profile-constraints}
Let \(s=\frac{\alpha(G[J])}{n}\) and \(\tau=1-r-s\).
Then $0\le s\le r$, $\tau\ge0$, and
\begin{equation}\label{eq:profile-constraints}
 z\le s\tau,\qquad w\le rs,\qquad
 z+w\le\frac r2,\qquad 2z+w\le r(1-r).
\end{equation}
\end{lemma}

\begin{proof}
The first inequality is Lemma~\ref{lem:alphaedges} applied to $G[J]$.
For every $v\in I$, the set $N(v)\cap J$ is independent in $G[J]$, so
$v$ has at most $sn$ neighbors in $J$; summing over $I$ gives $w\le rs$.  Finally,
$\Delta(G)\le\alpha(G)=rn$.  Summing degrees over all of $V(G)$ gives
$z+w\le r/2$, while summing over $J$ gives
$2z+w\le r(1-r)$.
\end{proof}

\paragraph{A sampling obstruction.}
More generally, one can choose prescribed proportions independently from
$I$ and $J$ while keeping the two sides equal.  The resulting family of
main-term edge-count vectors contains both canonical samples.
Even this full family of main-term vectors can miss the target, as the
next example shows. Appendix~\ref{app:full-sampling-obstruction} gives the
parametrization and exact minimization. The dichotomy itself uses only
the two costs in \eqref{eq:canonical-costs}.

\begin{example}
\label{ex:sampling-obstruction}
The exact feasible profile
\[
 r=\frac{10}{21},\qquad s=\tau=\frac{11}{42},\qquad
 z=\frac{121}{1764},\qquad w=\frac{11}{98}
\]
satisfies equality in the first and fourth constraints of
\eqref{eq:profile-constraints}.  Every main-term vector in the two-sided
sampling family has normalized $\ell_2$ cost at least
\[
 \sqrt{\frac{156269119}{39829708800}}>\frac1{16}.
\]
Appendix~\ref{app:full-sampling-obstruction} gives the exact minimizer.
Thus the entire main-term family misses even the limiting $1/16$ target;
compactness and continuity give the same obstruction for all
$p\in(2,\infty)$ sufficiently close to $2$.

The profile is realized asymptotically by an explicit five-class
triangle-free construction, yet those graphs admit bisections with only
$O(n)$ internal edges; see Remark~\ref{rem:obstruction-realization}.
\end{example}

\begin{remark}\label{rem:obstruction-realization}
For $n=42k$, take independent classes $A,D$ of size $9k$, classes $S,T$
of size $11k$, and $E=\{e_1,\dots,e_{2k}\}$.  Join $A$ completely to $S$
and $T$ completely to $D$; join $S$ completely to $T$ except for a
matching $s_it_i$ of size $2k$; and join $e_i$ to $s_i,t_i$. This
triangle-free graph has maximum independent set
$A\cup D\cup E$ of size $20k$: any independent set meeting both $S$ and
$T$ can use only matched exceptional pairs, while one meeting at most one
of these classes has size at most $20k$.  It asymptotically realizes the
profile in Example~\ref{ex:sampling-obstruction} and has a bisection
obtained from $A\cup T$ and $S\cup D$ by splitting $E$, with only $O(n)$
internal edges.
\end{remark}

The obstruction concerns upper bounds obtained from expected main-term
vectors; it gives no lower bound on individual sampling outcomes.
When both canonical costs exceed the target, the resulting constraints
force missing-edge structure from which a different bisection can be
constructed.

\paragraph{Defect orientation and balancing.}

Let \(S\) be a maximum independent set of \(G[J]\), and put
\(T=J\setminus S\). The missing edges between $S$ and $T$ bound the products of the two
neighbor counts for vertices in $I$. Placing each vertex of $I$ on the
side where it has fewer neighbors
converts this control into a bound on the internal-edge count.

\begin{lemma}
\label{lem:defect-orientation}
Let $I,S,T$ be pairwise disjoint, and assume that no triangle has one
vertex in each set.  Write
\[
 \mathcal M_{ST}=\{(s,t)\in S\times T:st\notin E(G)\},\qquad
 M_{ST}=|\mathcal M_{ST}|.
\]
For $v\in I$, put $a_v=|N(v)\cap S|$, $b_v=|N(v)\cap T|$ and
$m_v=\min\{a_v,b_v\}$.  Then
\[
 a_vb_v\le M_{ST}\quad(v\in I),\qquad
 \sum_{v\in I}m_v\le |I|\sqrt{M_{ST}}.
\]
If $I$ is independent and $P=\{v\in I:a_v\le b_v\}$, then
\[
 e(S\cup P)+e(T\cup(I\setminus P))
 =e(S)+e(T)+\sum_{v\in I}m_v.
\]
\end{lemma}

\begin{proof}
For each $v\in I$, every pair in $(N(v)\cap S)\times(N(v)\cap T)$ is
a missing $S$--$T$ edge. Hence
\(m_v\le\sqrt{a_vb_v}\le\sqrt{M_{ST}}\), and summing over $I$ gives the asserted bound.
If $I$ is independent, placing each $v$ on the side on which it has fewer
neighbors gives the stated identity.
\end{proof}

The next lemma quantifies the cost of restoring exact balance. Its defects
$\delta_1,\delta_2$ measure missing edges and degree loss, while $q$
records the shortfall from an independent half.

\begin{lemma}\label{lem:structural-bisection}
Assume $s\ge1/2-r$ and $s,\tau>0$.  Define
\[
 \delta_1=s\tau-z,\qquad
 \delta_2=r(1-r)-(2z+w),\qquad
 L=r\sqrt{2\delta_1},\qquad q=1-2r.
\]
Then $G$ has a bisection $V(G)=A\cup B$ such that
\[
 \frac{e(A)+e(B)}{n^2}\le \mathcal B(r,s,z,w),
\]
where
\begin{equation}\label{eq:def-B}
\mathcal B=\max\left\{
 \delta_1+L+s\left(\frac q2+\frac{\delta_2+L}{\tau}\right),
 \delta_1+L+\tau\left(\frac q2+\frac{\delta_2+L}{s}\right)
 \right\}.
\end{equation}
\end{lemma}

\begin{proof}
Let $S$ be a maximum independent set of $G[J]$ and put $T=J\setminus S$.
Thus $|S|=sn$ and $|T|=\tau n$.  Write
\[
 E_T=e(T),\qquad M_{ST}=|S||T|-e(S,T).
\]
The degree calculation behind Lemma~\ref{lem:alphaedges} gives
$M_{ST}\ge2E_T$, while
\[
 \delta_1n^2=M_{ST}-E_T.
\]
Therefore
\begin{equation}\label{eq:structural-first-defect}
 E_T\le\delta_1n^2,\qquad M_{ST}\le2\delta_1n^2.
\end{equation}

For $v\in I$, put $a_v=d_S(v)$ and $b_v=d_T(v)$. Lemma~\ref{lem:defect-orientation}
and \eqref{eq:structural-first-defect} give
\begin{equation}\label{eq:minority-degree}
 \sum_{v\in I}\min\{a_v,b_v\}\le rn\sqrt{M_{ST}}\le Ln^2.
\end{equation}

Let $P=\{v\in I:a_v\le b_v\}$.  Place $S\cup P$ on one side and
$T\cup(I\setminus P)$ on the other.  Before balancing, the total number of
internal edges is at most $(\delta_1+L)n^2$ by
\eqref{eq:structural-first-defect}--\eqref{eq:minority-degree}.

To restore balance, we compare the size of the oriented class with its required size. Put
\[
 p_0=\frac{|P|}{n},\qquad
 W_S=\frac{e(I,S)}{n^2},\qquad W_T=\frac{e(I,T)}{n^2},
 \qquad k=\frac12-s=r-s+\frac q2.
\]
Since $rn-d(v)\ge0$ for every $v\in J$, the defect $\delta_2n^2$ is the
sum of these degree deficits over $J$.  The deficit on $S$ is at most the
total deficit, and hence
\[
 e(I,S)+e(S,T)\ge(rs-\delta_2)n^2.
\]
Using $e(S,T)\le s\tau n^2$ gives $W_S\ge s(r-\tau)-\delta_2$.  Similarly,
\[
 e(I,T)+e(S,T)+2e(T)\ge(r\tau-\delta_2)n^2,
\]
whereas the maximum-independent-set degree bound in $G[J]$ gives
$e(S,T)+2e(T)\le s\tau n^2$.  Therefore
\begin{equation}\label{eq:orientation-lower}
 W_S\ge s(r-\tau)-\delta_2,\qquad
 W_T\ge \tau(r-s)-\delta_2.
\end{equation}

If $p_0<k$, then
$W_T\le \tau p_0+L$: vertices in $P$ have $T$-degree at most $\tau n$, and the
$T$-degrees of vertices in $I\setminus P$ are included in
\eqref{eq:minority-degree}.  Together with
\eqref{eq:orientation-lower}, this gives
\[
 k-p_0\le\frac q2+\frac{\delta_2+L}{\tau}.
\]
Move exactly $(k-p_0)n$ vertices from $I\setminus P$ to the $S$-side. This
number is an integer, and there are enough available vertices because
$k\le r$. Since $I$ is independent, each move increases the total
internal-edge count by at most $sn$. If $p_0>k$, the symmetric inequalities
$W_S\le L+s(r-p_0)$ and \eqref{eq:orientation-lower} give
\[
 p_0-k\le\frac q2+\frac{\delta_2+L}{s},
\]
and each required move costs at most $\tau n$. Both final classes have size
exactly $n/2$. These are precisely the two expressions in
\eqref{eq:def-B}, with no rounding or asymptotic error.
\end{proof}

\paragraph{Localization of the double-failure regime.}

For abstract feasible quadruples, we use the same formulas in
\eqref{eq:canonical-costs}, including at $r=1/2$.
The next lemma shows that if both $C_{\rm rnd}$ and $C_{\rm fix}$
exceed $c_*$, the structural defects lie in the range where orientation
and balancing succeed.

\begin{samepage}
\begin{lemma}
\label{lem:canonical-localization}
Let $(r,s,z,w)$ be a nonnegative quadruple satisfying
\eqref{eq:profile-constraints}, with $0<r\le1/2$, $0\le s\le r$, and
$\tau=1-r-s\ge0$.  Put
\[
 q=1-2r,\qquad \delta_1=s\tau-z,\qquad
 \delta_2=r(1-r)-(2z+w),
\]
\[
 \Delta=|s-\tau|,\qquad \omega=\frac{\Delta^2}{4}+\delta_1.
\]
If $C_{\rm rnd}>c_*$ and $C_{\rm fix}>c_*$, then
\[
 s\ge\frac12-r,\qquad s,\tau>0,
\]
and
\[
 q<\frac1{12},\qquad \omega<\frac1{1000},\qquad
 0\le\delta_2<h(q):=\frac1{80}-\frac q8-\frac{5q^2}{16}.
\]
\end{lemma}
\end{samepage}

\begin{proof}
Put $m=z+w$.  Assume both canonical costs exceed $c_*$.  Random-sample failure gives
\[
 m>2\sqrt2\,c_*>\frac{22}{125},
\]
where the last inequality follows by squaring $623\sqrt2>880$.
As $m\le r/2=(1-q)/4$, we obtain $q<37/125$.
If $s<q/2$, then $s<3/20$ and
\[
 m\le s\tau+rs=s(1-s)<\frac3{20}\frac{17}{20}
 =\frac{51}{400}<\frac{22}{125},
\]
a contradiction; here $x(1-x)$ is increasing on $[0,1/2]$.
Consequently $s\ge q/2=1/2-r$.  Also $s=0$ forces $z=w=0$,
contradicting random-sample failure.  If $\tau=0$, then $s\le r\le1/2$ forces
$s=r=1/2$, whence $q=z=0$ and $C_{\rm fix}=0$, another contradiction.
Thus $s,\tau>0$, as needed for the divisions in the structural bound.

Since \(s+\tau=(1+q)/2\), the mass identities are
\begin{equation}\label{eq:canonical-masses}
 z=\frac{(1+q)^2}{16}-\omega,\qquad
 m=\frac{3-2q-5q^2}{16}+\omega-\delta_2.
\end{equation}
For the fixed-$I$ sample, $\lambda=q/(1+q)$ and $\mu=1/(1+q)$,
so its coordinates are
\[
 X=\lambda m-\lambda\mu z,\qquad Y=\mu^2z.
\]
They are nonnegative at every feasible profile. At fixed $q,\omega$,
lowering $\delta_2$ to zero increases $X$ by $q\delta_2/(1+q)$ and
leaves $Y$ unchanged. When $\delta_2=0$, both coordinates are affine in
$\omega$, so $X^2+Y^2$ is convex in $\omega$. The endpoint checks in
Appendix~\ref{app:canonical-localization-proof} exclude $q\ge1/12$ and then
$\omega\ge1/1000$; the random-sample inequality then bounds $\delta_2$.
Only an algebraic upper envelope is used: the modified coordinates need not
be realized by another graph.
\end{proof}

\begin{lemma}[Sampling-or-structure dichotomy]
\label{lem:canonical-dichotomy}
Every nonnegative quadruple $(r,s,z,w)$ satisfying
\eqref{eq:profile-constraints}, with $0<r\le1/2$, $0\le s\le r$, and
$\tau=1-r-s\ge0$, satisfies at least one of the following alternatives:
\begin{enumerate}
\item $C_{\rm rnd}\le c_*$;
\item $C_{\rm fix}\le c_*$;
\item $s\ge1/2-r$, $s,\tau>0$, and $\mathcal B(r,s,z,w)<c_*$.
\end{enumerate}
Consequently, every triangle-free graph $G$ with $\alpha(G)<n/2$
satisfies $\Db_2(G)\le c_*n^2+O(n^{3/2})$.
\end{lemma}

\begin{proof}
Assume that the first two alternatives fail.  By
Lemma~\ref{lem:canonical-localization}, the structural domain holds and,
with
\[
 q=1-2r,\qquad \Delta=|s-\tau|,\qquad
 \omega=\frac{\Delta^2}{4}+\delta_1,
\]
we have
\[
 q<\frac1{12},\qquad \omega<\varepsilon:=\frac1{1000},\qquad
 0\le\delta_2<h(q).
\]
Write
\[
 b=1-r,\qquad s_+=\max\{s,\tau\},\qquad s_-=\min\{s,\tau\},
 \qquad R=\frac{s_+}{s_-},\qquad L=r\sqrt{2\delta_1}.
\]
The two terms in \eqref{eq:def-B} differ by
$(R-R^{-1})(L+\delta_2)+q(s_+-s_-)/2\ge0$; hence the worse orientation is
\begin{equation}\label{eq:canonical-worse-orientation}
 \mathcal B=\delta_1+(1+R)L+R\delta_2+\frac q2s_+.
\end{equation}
Set $\rho=2\sqrt\varepsilon=\sqrt{10}/50$.  Since $\omega<\varepsilon$,
we have $\Delta<\rho<b$ and therefore
\[
 R\le R_*(q):=\frac{b+\rho}{b-\rho},\qquad
 s_+\le\frac{b+\rho}{2}.
\]
The following identity couples the class imbalance with the missing-edge
defect:
\[
 \varepsilon(b-\Delta)^2
 -\left(\varepsilon-\frac{\Delta^2}{4}\right)(b^2-4\varepsilon)
 =\frac{(b\Delta-4\varepsilon)^2}{4}\ge0.
\]
Since $\delta_1<\varepsilon-\Delta^2/4$, this gives
\[
 (1+R)L<\frac{2br\sqrt{2\varepsilon}}
 {\sqrt{b^2-4\varepsilon}}.
\]
Substituting these estimates and $\delta_2<h(q)$ into
\eqref{eq:canonical-worse-orientation} yields
\begin{equation}\label{eq:canonical-envelope}
 \mathcal B<H(q):=\varepsilon
 +\frac{2br\sqrt{2\varepsilon}}{\sqrt{b^2-4\varepsilon}}
 +R_*(q)h(q)+\frac q4(b+\rho).
\end{equation}
Appendix~\ref{app:canonical-dichotomy-proof} proves that \(H\) is decreasing
on \(0\le q<1/12\) and that \(c_*-H(0)>17/1240000\).
Thus \(\mathcal B<c_*\), proving the third alternative.

For a graph, in either sampling alternative the corresponding canonical
procedure from the start of this section gives $\ell_2$ cost at most $c_*n^2+O(n^{3/2})$; the structural
alternative gives an exact bisection of $\ell_1$ cost below $c_*n^2$
by Lemma~\ref{lem:structural-bisection}.  Since
$\|\cdot\|_2\le\|\cdot\|_1$, the graph consequence follows.
\end{proof}

The following margins compare the structural bound with the model value
and supply the density thresholds for the central-contact range.

\begin{lemma}\label{lem:phase-rational-margins}
With \(c_*=623/10000\) as above,
\[
 \frac83<p_c<\frac{11}{4},\qquad
 c_*+10^{-5}<\gamma_{p_c}.
\]
For every $p_c\le p<\infty$, if \(\rho_p=\sigma_p(1-\sigma_p)\), then
\[
 \rho_p>\frac{481}{2500}>\frac{19}{100}.
\]
 At $p=4$ one has
 \[
  \frac27<\sigma_4<\frac5{17},\qquad
  \rho_4>\frac{10}{49},\qquad
  \gamma_4>\frac{35}{578}>\frac3{50}+10^{-5}.
 \]
\end{lemma}

\begin{proof}
The exact rational comparisons are collected in
Appendix~\ref{app:rational-margins}.
\end{proof}
\begin{theorem}\label{thm:interior-range}
For every fixed $2<p<p_c$, every triangle-free graph $G$ of even order $n$
with $\alpha(G)<n/2$ satisfies
\begin{equation}\label{eq:uniform-alpha-gap}
 \Db_p(G)\le c_*n^2+O(n^{3/2})
 \le(\gamma_p-10^{-5})n^2+O(n^{3/2}).
\end{equation}
\end{theorem}

\begin{proof}
Form the maximum-independent-set profile from Section~\ref{sec:unified} and
apply Lemma~\ref{lem:canonical-dichotomy}. Its graph consequence and norm
monotonicity give
\[
 \|(e(A),e(B))\|_p\le\|(e(A),e(B))\|_2
 \le c_*n^2+O(n^{3/2}).
\]
Since $\gamma_p>\gamma_{p_c}>c_*+10^{-5}$ by
Corollary~\ref{cor:monotone} and Lemma~\ref{lem:phase-rational-margins}, this
proves \eqref{eq:uniform-alpha-gap}. For sufficiently large $n$, the fixed gap absorbs the sampling error.
\end{proof}

\subsection{Central contact: the density reduction}

For $p_c\le p<\infty$ the optimal contact is central, and one
independent-set density argument handles this range.
The split at $p=4$ below only changes the comparison norm.

\begin{lemma}\label{lem:central-gap}
Let $p_c\le p<\infty$, and let $G$ be triangle-free with
$\alpha(G)<n/2$. Then
\[
 \Db_p(G)\le(\gamma_p-10^{-5})n^2+O(n^{3/2}).
\]
\end{lemma}

\begin{proof}
Write $\rho_p=\sigma_p(1-\sigma_p)$. The central-contact formula gives
$\gamma_p=2^{1/p}\rho_p/4$.
First suppose $p_c\le p\le4$ and put $\rho_0=19/100$.  By
Lemma~\ref{lem:phase-rational-margins},
\[
 \gamma_p-\frac{2^{1/p}\rho_0}{4}
 =\frac{2^{1/p}(\rho_p-\rho_0)}4
 >\frac3{5000}.
\]
If $e(G)\le\rho_0n^2$, Lemma~\ref{lem:balanced-sampling} gives the required gap.
Otherwise a maximum independent set has size $rn$ with
$r>19/50$.  Put $q=1-2r\in(0,6/25)$.  Lemma~\ref{lem:fixedI}, with
$\rho=\rho_0$, gives normalized coordinates bounded, up to a common
$O(n^{-1/2})$ error, by
\[
 A=\frac{q(1-q)}{4(1+q)},\qquad
 B_0=\frac{6/25-q^2}{4(1+q)^2}.
\]
For $c_0=3/50$ these satisfy the coordinatewise estimate
\[
 (A,B_0)\le\frac1{(1+q)^2}\left(\frac q4,c_0\right).
\]
Indeed the first coordinate reduces to $1-q^2\le1$, and the second to
$6/100=c_0$.  Since $q<1/4$,
\[
 (1+q)^4\ge1+22q^2,
 \qquad \frac1{16}<22c_0^2.
\]
The first inequality follows from $4q\ge16q^2$ after expanding the left
side.  Squaring the coordinatewise bound now gives
$\|(A,B_0)\|_2<c_0$.  Hence
\[
 \|(A,B_0)\|_p\le\|(A,B_0)\|_2<\frac3{50}
 <\gamma_4-10^{-5}\le\gamma_p-10^{-5}.
\]

Now let $4\le p<\infty$ and put $\rho_0=1/5$.
Corollary~\ref{cor:monotone} gives $\sigma_p\ge\sigma_4$.
Since $s(1-s)$ is increasing on $[0,1/2]$,
$\rho_p\ge\rho_4>10/49$.  Thus the random branch
$e(G)\le\rho_0n^2$ lies below $\gamma_p$ by more than
\[
 \frac14\left(\frac{10}{49}-\frac15\right)=\frac1{980}.
\]
In the other branch $r>2/5$, and Lemma~\ref{lem:fixedI}, with
$\rho=1/5$, gives
\[
 A(r)\le\frac1{30},\qquad B(r)\le\frac1{20}.
\]
The first bound follows by differentiating $A(r)$ on $[2/5,1/2]$; the
second is equivalent to
$(3r-2)(2r-1)\ge0$.  Corollaries~\ref{cor:monotone} and~\ref{cor:phase-explicit} also give
$\gamma_p\ge\gamma_\infty=1/18$.  Therefore, for every
$p\ge4$,
\[
 \|(A(r),B(r))\|_p\le\|(A(r),B(r))\|_4
 \le\frac{\sqrt[4]{97}}{60}
 <\frac4{75}<\frac1{18}-10^{-5}\le\gamma_p-10^{-5}.
\]
Here $\sqrt[4]{97}<16/5$ follows by taking fourth powers. In both ranges,
$\|v\|_p\le\|v\|_1$ bounds the norm of the normalized sampling error by
$O(n^{-1/2})$. Adding these errors gives the stated bound.
\end{proof}

\begin{corollary}
\label{cor:alphagap}
For every fixed $2<p<\infty$ there is $n_p$ such that every
triangle-free $G$ of even order $n\ge n_p$ with $\alpha(G)<n/2$ satisfies
\[
 \frac{\Db_p(G)}{n^2}\le\gamma_p-\frac{10^{-5}}2.
\]
\end{corollary}

\begin{proof}
For $2<p<p_c$, use the strict estimate
\eqref{eq:uniform-alpha-gap}; for $p_c\le p<\infty$, use
Lemma~\ref{lem:central-gap}. Both estimates have gap $10^{-5}$ and
normalized error $O_p(n^{-1/2})$. Increase $n_p$ until this error is at
most $10^{-5}/2$.
\end{proof}

\section{Quantitative stability}\label{sec:stability}

For a near-extremal graph, Corollary~\ref{cor:alphagap} rules out
$\alpha(G)<n/2$ and hence forces an independent half. The profile corner
locates the smaller class, while
sampling controls the edge-count slack. Triangle-freeness then turns
these slack bounds into an edit-distance estimate.

Fix $2<p<\infty$ and write $\omega_n=n^{-1}$.
The $O_p(\omega_n)$ scalar error is supplied by
Lemma~\ref{lem:squared-norm-sampling} under near-extremality, rather than
a coordinatewise sampling bound. The coarser sampling errors have already
been absorbed into the fixed small-independence gap.

\begin{lemma}
\label{lem:independenthalf-stability}
Let $I$ be an independent set of size $n/2$, put $J=V(G)\setminus I$, and
let $S$ be a maximum independent set of $G[J]$.  Write
\[
 |S|=yn,\quad b=\frac{e(J)}{n^2},\quad w=\frac{e(I,J)}{n^2},\quad
 \Delta_J=y(1/2-y)-b,\quad \Delta_{IJ}=y/2-w.
\]
There is a constant $A_p$ such that, if
$\Db_p(G)\ge(\gamma_p-\delta)n^2$ with $0\le\delta\le\gamma_p/2$, then
\[
 |y-\sigma_p|\le\frac{\delta+A_p\omega_n}{c_p^{\rm prof}}.
\]
Moreover, put
\[
 \varrho_p=\frac{\sigma_p-1/4}{2},\qquad
 c_p^{\rm slack}=
 \frac{\theta_p(1-\theta_p)}p\,2^{-(p-1)/p},
 \qquad \theta_p=\frac{1-t_p}{2}.
\]
Whenever $|y-\sigma_p|\le\varrho_p$,
\[
 \Delta_J+\Delta_{IJ}\le
 \frac{\delta+A_p\omega_n+|y-\sigma_p|}{c_p^{\rm slack}}.
\]
\end{lemma}

\begin{proof}
Lemma~\ref{lem:alphaedges} gives $\Delta_J\ge0$, and the independence of
$N(v)\cap J$ gives $\Delta_{IJ}\ge0$. Both conclusions will follow from
a sampling estimate valid for every $\theta\in[0,1]$. Choose an integer
$k=\theta n/2+O(1)$ and independent uniform $k$-sets
$P\subseteq I$ and $U\subseteq J$. The sets
$(I\setminus P)\cup U$ and $P\cup(J\setminus U)$ form a bisection.
Its normalized mean vector is $z_\theta+O(n^{-1})$, where
\[
 z_\theta=
 \bigl(\theta^2b+\theta(1-\theta)w,\,
 (1-\theta)^2b+\theta(1-\theta)w\bigr)\ge0.
\]
Put $a=\gamma_p-\delta\ge\gamma_p/2$. By
Lemma~\ref{lem:squared-norm-sampling},
\[
 a^2\le\|z_\theta\|_p^2+C_p/n.
\]
If $\|z_\theta\|_p<a$, factor the difference of squares and use
$a+\|z_\theta\|_p\ge\gamma_p/2$; otherwise the following inequality
is immediate. Choosing $A_p\ge2C_p/\gamma_p$, we obtain, uniformly in
$\theta$,
\begin{equation}\label{eq:independent-half-sampling}
 \gamma_p-\delta\le\|z_\theta\|_p+A_p\omega_n.
\end{equation}

Choose a minimizing parameter $\theta_y$ for $F_p(y)$. Since
$z_{\theta_y}\le(\mathcal C_1(y,\theta_y),\mathcal C_2(y,\theta_y))$ coordinatewise,
\eqref{eq:independent-half-sampling} gives
$\gamma_p-\delta\le F_p(y)+A_p\omega_n$.
The linear profile bound~\eqref{eq:linear-profile-corner} proves the
first assertion.

For the slack bound, take $\theta=\theta_p$ in \eqref{eq:independent-half-sampling}, and put
\[
 x=(\mathcal C_1(y,\theta_p),\mathcal C_2(y,\theta_p)),\qquad z=z_{\theta_p}=x-d,
\]
where
\[
 d=\bigl(\theta_p^2\Delta_J+\theta_p(1-\theta_p)\Delta_{IJ},\,
 (1-\theta_p)^2\Delta_J+\theta_p(1-\theta_p)\Delta_{IJ}\bigr).
\]
Thus $0\le z\le x$ coordinatewise. Since $\theta_p\le1/2$ and $y\le1/2$,
we have $x_2\ge x_1$ and $x_2/\|x\|_p\ge2^{-1/p}$. Also
\[
 d_2\ge\theta_p(1-\theta_p)(\Delta_J+\Delta_{IJ}).
\]
Concavity of the $p$th root and
$x_i^p-(x_i-d_i)^p\ge x_i^{p-1}d_i$ give
\[
 \|x\|_p-\|z\|_p
 \ge\frac{x_1^{p-1}d_1+x_2^{p-1}d_2}{p\|x\|_p^{p-1}}
 \ge c_p^{\rm slack}(\Delta_J+\Delta_{IJ}).
\]
Each $\mathcal C_i(\cdot,\theta_p)$ is $1/2$-Lipschitz, and the contact equation
gives $\|x(\sigma_p)\|_p=\gamma_p$. Hence
$\|x(y)\|_p\le\gamma_p+|y-\sigma_p|$, since $2^{1/p-1}\le1$.
Combining this with \eqref{eq:independent-half-sampling} at $\theta_p$ proves
\[
 c_p^{\rm slack}(\Delta_J+\Delta_{IJ})
 \le\delta+A_p\omega_n+|y-\sigma_p|.
\]
\end{proof}

\begin{lemma}[Homogeneous defect bound]\label{lem:slack-edit}
In the notation of Lemma~\ref{lem:independenthalf-stability}, assume
\(|y-\sigma_p|\le(\sigma_p-1/4)/2\). Put \(T=J\setminus S\),
\(X=S\), \(Y=I\cup T\), and
\[
 h=e(Y),\qquad m=|X||Y|-e(X,Y),\qquad Q=m-h,\qquad D=m+h.
\]
Then \(X\) is independent and
\begin{equation*}
 Q=(\Delta_J+\Delta_{IJ})n^2\ge0,\qquad
 D\le E_pQ,\qquad E_p=3+\frac4{\sigma_p-1/4}.
\end{equation*}
In particular, \(Q=0\) implies \(G=K_{X,Y}\).
\end{lemma}

\begin{proof}
Let $M_{ST}=|S||T|-e(S,T)$.  The maximum-independent-set degree bound
gives $e(S,T)+2e(T)\le|S||T|$, whereas
$\Delta_Jn^2=M_{ST}-e(T)$.  Hence
\begin{equation}\label{eq:st-slack-bounds}
 e(T)\le\Delta_Jn^2,\qquad M_{ST}\le2\Delta_Jn^2.
\end{equation}
For $v\in I$, put $a_v=|N(v)\cap S|$, $b_v=|N(v)\cap T|$, and
$d_v=|S|-a_v-b_v\ge0$; the inequality follows because
$N(v)\cap J$ is independent and $S$ is maximum.  Then
\[
 \sum_{v\in I}d_v=\Delta_{IJ}n^2,
\]
and Lemma~\ref{lem:defect-orientation} gives $a_vb_v\le M_{ST}$.

Write $g=|S|-|T|=(2y-1/2)n\ge(\sigma_p-1/4)n>0$.  Pointwise,
\[
 a_v+d_v=|S|-b_v\ge g,
\]
and therefore
\[
 g b_v\le(a_v+d_v)b_v
 \le M_{ST}+|T|d_v.
\]
Summing this inequality over $I$ and using
\eqref{eq:st-slack-bounds} gives the direct estimate
\[
 \frac{e(I,T)}{n^2}\le
 \frac{\Delta_J+(1/2-y)\Delta_{IJ}}{2y-1/2}
 \le\frac{\Delta_J+\Delta_{IJ}/2}{\sigma_p-1/4}.
\]
Exact counting gives
\[
 h=e(T)+e(I,T),\qquad m=M_{ST}+\Delta_{IJ}n^2+e(I,T).
\]
Since \(\Delta_Jn^2=M_{ST}-e(T)\), this proves the identity for \(Q\).
The edit defect \(D\) is
\(e(T)+M_{ST}+\Delta_{IJ}n^2+2e(I,T)\). The preceding bounds imply
\[
 \frac{e(X)+e(Y)+|X||Y|-e(X,Y)}{n^2}
 \le\left(3+\frac4{\sigma_p-1/4}\right)
 (\Delta_J+\Delta_{IJ}),
\]
which proves the homogeneous bound. Since \(X=S\) is independent,
\(D=0\) is equivalent to the graph being complete bipartite on \((X,Y)\).
\end{proof}

\begin{proof}[Proof of Theorem~\ref{thm:intro-stability}]
Put $\eta_0=10^{-5}$, fix the preceding constants, and enlarge $A_p$ if
necessary.  Choose
\[
 0<\delta_p\le\frac14
 \min\{\eta_0,\gamma_p,c_p^{\rm prof}\varrho_p\},\qquad
 A_p\omega_n\le\frac14c_p^{\rm prof}\varrho_p
 \quad(n\ge n_p),
\]
where $n_p$ also satisfies Corollary~\ref{cor:alphagap}.
That corollary rules out $\alpha(G)<n/2$, since otherwise
$\delta\ge\eta_0/2$, contrary to $\delta\le\eta_0/4$.
Take an independent half $I$.
The first assertion of Lemma~\ref{lem:independenthalf-stability} gives
$|y-\sigma_p|\le\varrho_p/2$ by the choice of $\delta_p$ and $n_p$.
Its second assertion and
Lemma~\ref{lem:slack-edit} therefore apply and give
\[
 |y-\sigma_p|\le\frac{\delta+A_p\omega_n}{c_p^{\rm prof}},\qquad
 \Delta_J+\Delta_{IJ}\le
 \frac{(1+1/c_p^{\rm prof})(\delta+A_p\omega_n)}{c_p^{\rm slack}}.
\]
The first inequality controls the class-size term in
Theorem~\ref{thm:intro-stability}. The second, together with
Lemma~\ref{lem:slack-edit}, controls the normalized edit defect of the
partition $(S,I\cup T)$. Absorbing the fixed $p$-dependent coefficients
into $C_p$ proves the first assertion. To fix the class size, move at most
$|y-\sigma_p|n+1$ vertices to impose the class size
$\lfloor\sigma_pn\rfloor$.  Each move changes at most $n-1$ target
pairs, so its additional normalized cost is at most
$|y-\sigma_p|+1/n=O_p(\delta+\omega_n)$.  This proves the second assertion.
\end{proof}

For $p=2$, the right-branch argument in
Proposition~\ref{prop:bipartite-profile} remains valid,
since $\kappa_2=1/2$. It follows that
\[
F_2(1/4+\varepsilon)=\frac1{16}-\varepsilon^2
\qquad(0\le\varepsilon\le1/4).
\]

\section{Exact extremizers}\label{sec:exactification}

Fix $2<p<\infty$, write $N=n/2$, and use
\eqref{eq:discrete-model}.
We now turn quantitative stability into an exact
classification. For a global extremizer, the preceding
estimates give a nonnegative integer defect $Q=O_p(n)$
and an edit bound proportional to $Q$.
Sampling at an integer-minimizing quota then shows
that every $Q>0$ incurs a strict loss, including
when the minimum is attained at several quotas
or at an endpoint.

\begin{lemma}\label{lem:exact-defects}
Let \(G\) maximize \(\Db_p\) among triangle-free graphs of even order \(n\).
For sufficiently large \(n\), an independent half exists. Choose one, \(I\), put
\(J=V(G)\setminus I\), and let \(S\) be any maximum independent set of \(G[J]\).
Set \(T=J\setminus S\), \(X=S\), \(Y=I\cup T\), \(a=|S|\),
\(h=e(Y)\), and \(m=a(n-a)-e(X,Y)\). Then
\[
 a=\sigma_p n+O_p(1),\qquad
 Q:=m-h\ge0,\qquad Q\le K_p n,\qquad D:=m+h\le E_pQ,
 \quad E_p=3+\frac4{\sigma_p-1/4},
\]
for a constant \(K_p\) depending only on \(p\).
In particular, \(Q=0\) implies \(G=K_{X,Y}\).
\end{lemma}
\begin{proof}
Put \(L_0=2^{1/p-1}\). Proposition~\ref{prop:bipartite-profile} gives
\(\Db_p(G)\ge\gamma_p n^2-L_0n\), so Corollary~\ref{cor:alphagap}
forces an independent half for large \(n\). Apply
Lemma~\ref{lem:independenthalf-stability} with \(\delta=L_0/n\).
For large enough \(n\), its first estimate places \(y\) in the neighborhood
required for the slack bound. Together, the two estimates give
\[
 |a-\sigma_pn|\le C_p^{\rm loc}:=\frac{L_0+A_p}{c_p^{\rm prof}},\qquad
 (\Delta_J+\Delta_{IJ})n^2\le K_p n,
\]
where one may take
\[
 K_p=\frac{(1+1/c_p^{\rm prof})(L_0+A_p)}{c_p^{\rm slack}}.
\]
Lemma~\ref{lem:slack-edit} applies to the same sets \(X=S\),
\(Y=I\cup T\). Its exact identity and homogeneous bound give all the
remaining assertions, including the conclusion when \(Q=0\).
\end{proof}

\begin{lemma}[Integer quotas and weighted defect loss]\label{lem:exact-loss}
In the notation of Lemma~\ref{lem:exact-defects}, let \(k\) be any integer
minimizer defining \(H_{p,n}(a)\), and set
\[
 C=(k(N-k),(a-k)(N-a+k)),\quad R=\|C\|_p,\quad
 \theta=k/a,\quad b=n-a,\quad \eta=(N-k)/b.
\]
Choose independently uniform sets \(P\subseteq X\), \(U\subseteq Y\)
of sizes \(k,N-k\), and let \(Z\) be the internal-edge-count vector of
the bisection \((P\cup U,V(G)\setminus(P\cup U))\).
For sufficiently large \(n\), put
\(g=\left.\nabla_x\|x\|_p\right|_{x=C}\). Then
\[
 \|Z-C\|_2\le D\quad\hbox{in every outcome},\qquad
 g\cdot(C-\mathbb EZ)\ge Q/48.
\]
\end{lemma}
\begin{proof}
Every outcome has halves of size \(N\). For the complete bipartite
reference the vector is always \(C\). Each of the \(D\) defect edges
changes at most one coordinate by one, so
\(\|Z-C\|_2\le\|Z-C\|_1\le D\).
Exact fixed-size sampling gives
\begin{equation}\label{eq:exact-means}
\begin{aligned}
 C_1-\mathbb EZ_1&=m\theta\eta
       -h\left(\eta^2-\frac{\eta(1-\eta)}{b-1}\right),\\
 C_2-\mathbb EZ_2&=m(1-\theta)(1-\eta)
       -h\left((1-\eta)^2-\frac{\eta(1-\eta)}{b-1}\right).
\end{aligned}
\end{equation}
No coordinatewise comparison of \(C\) and \(\mathbb EZ\) is required.

By exchanging the halves, assume \(k\le a/2\). For large \(n\),
\(s=a/n\in[1/5,2/5]\). Using the coordinate maps from Section~\ref{sec:model}, set
\[
 c(s,v)=(\mathcal C_1(s,v),\mathcal C_2(s,v)),\qquad f_s(v)=\|c(s,v)\|_p.
\]
Then \(C=n^2c(s,\theta)\).
On this rectangle, \(c_1+c_2\ge s(1/2-s)\ge1/25\), so
\(f_s\ge1/50\). For fixed \(p>2\), this is a \(C^2\) function even
at \(v=0,1\): \(|x|^p\) is \(C^2\), and the vector never vanishes.
On this rectangle we have the uniform bound
\[
 |f_s''(v)|\le L_p:=4(p-1)+1.
\]
Indeed, the Hessian of the norm at a nonzero nonnegative vector is
positive semidefinite and at most \((p-1)I_2/\|c\|_p\);
each \(|\partial_v c_i|\le1/5\), and each
\(|\partial_v^2c_i|\le8/25\), while \(0\le g_i\le1\).
The assertion follows, including axes by continuity.

If \(1\le k\le a-1\), both neighbors are feasible in the original
integer domain, regardless of whether \(a\) is odd or \(k\) is central.
Taylor's inequality at the step \(1/a\) and the optimality of the
chosen quota therefore give
\[
 |f_s'(\theta)|\le L_p/(2a)\le 5L_p/(2n).
\]
The same bound holds for each integer minimizer, including tied minima.
The gradient is homogeneous of degree zero, so
\(g=\left.\nabla_x\|x\|_p\right|_{x=c(s,\theta)}\). Set \(d=\eta-\theta\),
\(G_0=g_1+g_2\), and
\[
 A_0=g_1\theta\eta+g_2(1-\theta)(1-\eta),\qquad
 H_0=g_2(1-\eta)-g_1\eta .
\]
Direct differentiation, using \((1-s)\eta=1/2-s\theta\), gives
\begin{equation}\label{eq:quota-gradient-identity}
 (1-2s)H_0=s\,d\,G_0-\frac{f_s'(\theta)}s.
\end{equation}
Since \(m=Q+h\), the mean identities \eqref{eq:exact-means}
imply the exact identity
\[
 g\cdot(C-\mathbb EZ)
 =A_0Q+h\left(dH_0+\frac{\eta(1-\eta)G_0}{b-1}\right).
\]
For interior quotas, \eqref{eq:quota-gradient-identity} gives
\[
 dH_0=\frac{s\,d^2G_0}{1-2s}
             -\frac{d\,f_s'(\theta)}{s(1-2s)}
 \ge-\frac{125L_p}{2n},
\]
using \(|d|\le1\), \(s\ge1/5\), and \(1-2s\ge1/5\).
Also \(C_2\ge C_1\), hence \(g_2\ge2^{-(p-1)/p}\ge1/2\).
As \(1-\theta\ge1/2\) and \(1-\eta\ge1/6\),
\(A_0\ge1/24\). The finite-population correction is nonnegative.
Using \(h\le D\le E_pQ\) yields
\[
 g\cdot(C-\mathbb EZ)\ge
 \left(\frac1{24}-\frac{125E_pL_p}{2n}\right)Q\ge Q/48
\]
when \(n\ge3000E_pL_p\).

If \(k=0\), then \(g=(0,1)\). The second exact mean formula instead gives
\[
 g\cdot(C-\mathbb EZ)
 =(1-\eta)Q+
 h\left(\eta(1-\eta)+\frac{\eta(1-\eta)}{b-1}\right)
 \ge Q/6 .
\]
The case \(k=a\) is symmetric. Hence both endpoints satisfy the asserted
lower bound without a stationarity estimate.
\end{proof}

\begin{proof}[Proof of Theorem~\ref{thm:exactification}]
Let \(G\) be a global extremizer. Lemmas~\ref{lem:exact-defects}
and~\ref{lem:exact-loss} apply for sufficiently large even \(n\).
Use the deterministic-reference form \eqref{eq:reference-sampling}
of Lemma~\ref{lem:squared-norm-sampling} with \(x=C\).
Since \(\nabla(\|\cdot\|_p^2)(C)=2R g\), its linear term is controlled by
Lemma~\ref{lem:exact-loss}.
Taking expectations and using \(R\ge n^2/50\) gives
\[
 \mathbb E\|Z\|_p^2
 \le R^2-\frac{RQ}{24}+(p-1)E_p^2Q^2
 \le R^2-\frac{n^2Q}{1200}+(p-1)E_p^2Q^2.
\]
For \(0<Q\le K_pn\), the right side is strictly less than \(R^2\)
once \(n>1200(p-1)E_p^2K_p\). This is a strict comparison for
every positive defect, including \(Q=1\). Since every sampled
outcome is a bisection,
\[
 \Db_p(G)^2\le\mathbb E\|Z\|_p^2<R^2.
\]
But \(R=H_{p,n}(a)\le B_p(n)\le\Db_p(G)\), a contradiction.
Thus \(Q=0\), and Lemma~\ref{lem:exact-defects} gives \(G=K_{a,n-a}\).
Its class size belongs to \(\mathcal A_p(n)\); conversely every listed
complete bipartite graph attains \(B_p(n)\). This establishes the value and all equality cases.

All thresholds above depend only on fixed \(p\). Choose \(n_p\) large
enough for the gap and slack inputs, \(L_0/n\le\gamma_p/2\),
\(C_p^{\rm loc}/n\le\min\{1/20,(\sigma_p-1/4)/2\}\), and the two explicit
thresholds. No uniformity as \(p\downarrow2\) or \(p\to\infty\)
is asserted. The argument applies at \(p_c\) as well: no nondegenerate
critical point or uniqueness of an integer minimizer has been used.
\end{proof}

\begin{corollary}\label{cor:integer-window}
Put \(L_0=2^{1/p-1}\) and \(w_p^{\rm int}=2L_0/c_p^{\rm prof}\).
For every even \(n\), the full set \(\mathcal A_p(n)\) is obtained by
maximizing \(H_{p,n}\) over
\[
 \mathcal W_p(n)=\{a\in\mathbb Z:0\le a\le n/2,\ |a-\sigma_p n|\le w_p^{\rm int}\}.
\]
For sufficiently large \(n\), these and only these maximizing class sizes
give the extremal graphs in Theorem~\ref{thm:exactification}.
\end{corollary}
\begin{proof}
For any integer \(a\),
every integer quota is feasible in the continuous profile, while rounding
a continuous minimizing value of \(x=s\theta\) to a feasible multiple of
\(1/n\) changes each coordinate by at most \(1/(2n)\).
Thus the bounds
\[
 n^2F_p(a/n)\le H_{p,n}(a)\le n^2F_p(a/n)+L_0n
\]
hold, because the coordinate maps are \(1/2\)-Lipschitz in \(x\).
If \(a\) maximizes \(H_{p,n}\), the rounded lower construction and the
linear profile corner imply
\[
 c_p^{\rm prof}|a-\sigma_p n|\,n
 \le (\gamma_p-F_p(a/n))n^2\le2L_0n.
\]
Hence every maximizing class size lies in \(\mathcal W_p(n)\),
which proves the complete discrete parameterization. Within that window,
one compares the finitely many exact values \(H_{p,n}(a)\), retaining
all equal maxima. The inner minimum can equivalently be taken over
\(0\le k\le\lfloor a/2\rfloor\) by symmetry.

\end{proof}

\begin{corollary}\label{cor:finite-asymptotic}
For every fixed \(2<p<\infty\), both estimates in
\eqref{eq:intro-asymptotic} hold as \(n\) tends to infinity through even integers.
\end{corollary}
\begin{proof}
Proposition~\ref{prop:bipartite-profile} bounds the integer complete
bipartite optimum by \(\gamma_p n^2+O_p(n)\) and supplies the rounded
construction. Apply Theorem~\ref{thm:exactification}.
\end{proof}

\appendix
\section{Analytic details for the support problem}\label{app:analytic}
\setcounter{equation}{0}
\renewcommand{\theequation}{\thesection.\arabic{equation}}
\renewcommand{\theHequation}{\thesection.\arabic{equation}}

\subsection{Curvature sign pattern}\label{app:curvature-sign}

The sign pattern in Lemma~\ref{lem:curvature} follows from the derivative
chain below. Termwise differentiation gives
\[
 Q_\alpha'(\xi)=\alpha\xi^{\alpha-2}R_\alpha(\xi),
\]
where
\[
\begin{aligned}
 R_\alpha(\xi)={}&-3\xi^{2\alpha+1}
 +(2\alpha+1)\xi-2(\alpha-1)\\
 &+2(2\alpha+1)\xi^{\alpha+2}
 -2(2\alpha+1)\xi^{\alpha+1}.
\end{aligned}
\]
Differentiating twice more gives
\[
 R_\alpha'=(2\alpha+1)K_\alpha,\qquad
 K_\alpha'=2\alpha\xi^{\alpha-1}H_\alpha,
\]
where
\[
 K_\alpha(\xi)=1-3\xi^{2\alpha}
 +2(\alpha+2)\xi^{\alpha+1}-2(\alpha+1)\xi^\alpha
\]
and
\[
 H_\alpha(\xi)=-(\alpha+1)
 +\frac{(\alpha+1)(\alpha+2)}{\alpha}\xi-3\xi^\alpha.
\]
The endpoint values are
\[
\begin{array}{c|cc}
 &0&1\\ \hline
 Q_\alpha&1&0\\
 R_\alpha&-2(\alpha-1)&0\\
 K_\alpha&1&0\\
 H_\alpha&-(\alpha+1)&(2-\alpha)/\alpha.
\end{array}
\]
For $1<\alpha<2$, moreover,
\[
 Q_\alpha'(1)=Q_\alpha''(1)=0,\qquad
 Q_\alpha'''(1)=2\alpha(2\alpha+1)(2-\alpha)>0,
\]
so, when $1<\alpha<2$, $\xi=1$ is an endpoint zero of exact multiplicity
three.

For $1<\alpha<2$ and $0\le\xi\le1$,
\[
 H_\alpha'(\xi)
 =\frac{(\alpha+1)(\alpha+2)}{\alpha}
  -3\alpha\xi^{\alpha-1}
 \ge\frac{(2\alpha+1)(2-\alpha)}{\alpha}>0.
\]
Thus $H_\alpha$ has a unique zero $\xi_H\in(0,1)$ and changes sign
from negative to positive.  It follows that $K_\alpha$ first strictly
decreases and then strictly increases to $K_\alpha(1)=0$. Hence it has a
unique zero $\xi_K\in(0,\xi_H)$ and changes sign from positive to negative.
Next $R_\alpha$ first strictly increases and then strictly decreases to
$R_\alpha(1)=0$, so it has a unique zero
$\xi_R\in(0,\xi_K)$ and changes sign from negative to positive.
Finally $Q_\alpha$ first strictly decreases and then strictly increases to
$Q_\alpha(1)=0$, so it has a unique zero
$\xi_Q\in(0,\xi_R)$ and changes sign from positive to negative. Therefore
\[
 0<\xi_Q<\xi_R<\xi_K<\xi_H<1.
\]

Now let $\alpha\ge2$. The weighted arithmetic--geometric mean
inequality, with weights $(\alpha-1)/\alpha$ and $1/\alpha$, gives
\[
 \alpha+1+3\xi^\alpha
 \ge\alpha\left(\frac{\alpha+1}{\alpha-1}\right)^{(\alpha-1)/\alpha}
 3^{1/\alpha}\xi.
\]
Thus $H_\alpha(\xi)\le0$ follows if the coefficient on the right is at
least $(\alpha+1)(\alpha+2)/\alpha$. Raising these positive coefficients
to the power $\alpha$ shows that this is equivalent to
\[
 T_\alpha:=
 \frac{(\alpha+1)(\alpha+2)^\alpha(\alpha-1)^{\alpha-1}}
 {\alpha^{2\alpha}}\le3.
\]
Factoring the powers of $\alpha$ gives
\[
 T_\alpha=\frac{(\alpha+1)(\alpha+2)}{\alpha^2}
 \left(1+\frac{\alpha-2}{\alpha^2}\right)^{\alpha-1}.
\]
Since $\log(1+x)\le x$ for $x\ge0$, we have the explicit chain
\[
\begin{aligned}
 \log T_\alpha
 &=\log\frac{(\alpha+1)(\alpha+2)}{\alpha^2}
    +(\alpha-1)\log\frac{(\alpha+2)(\alpha-1)}{\alpha^2}\\
  &=\log\frac{(\alpha+1)(\alpha+2)}{\alpha^2}
    +(\alpha-1)\log\left(1+\frac{\alpha-2}{\alpha^2}\right)\\
 &\le\log\frac{(\alpha+1)(\alpha+2)}{\alpha^2}
   +\frac{(\alpha-1)(\alpha-2)}{\alpha^2}\\
 &\le\log3.
\end{aligned}
\]
For the last inequality, the derivative of the first two terms on the right is
\[
 \frac{\alpha^2-6\alpha-8}
 {\alpha^3(\alpha+1)(\alpha+2)}.
\]
Its numerator changes sign exactly once, from negative to positive.
Hence that sum first decreases and then increases, so its supremum on
$[2,\infty)$ is the larger of its value $\log3$ at $2$
and its limit $1$. Since $\log3>1$, this proves
$T_\alpha\le3$, and hence $H_\alpha\le0$. Since
$K_\alpha(1)=R_\alpha(1)=Q_\alpha(1)=0$, the derivative chain gives
$K_\alpha\ge0$, $R_\alpha\le0$ and $Q_\alpha\ge0$, as required.

\subsection{Sign at the center}\label{app:center-sign}

Theorem~\ref{thm:contact-phase} uses the monotonicity and endpoint signs of
$f(p)=2a_p(a_p-1)-(p-1)$, where $a_p=2^{1-1/p}$. Since $a_p'=a_p\log2/p^2$,
$f'(p)=2a_p(2a_p-1)\log2/p^2-1$.
Set $G(p)=a_p(2a_p-1)/p^2$.  For $2<p<\infty$,
\[
 \frac{G'(p)}{G(p)}
 =\frac{\log2}{p^2}\frac{4a_p-1}{2a_p-1}-\frac2p
 <\frac3{p^2}-\frac2p<0.
\]
Moreover,
\[
 G(2)=1-\frac{\sqrt2}{4}<\frac23,\qquad \log2<\frac34.
\]
Here $\sqrt2>4/3$, while the first four positive terms of the exponential
series give $e^{3/4}>2$.  Hence $2\log2\,G(p)<1$ and $f'(p)<0$ on
$(2,\infty)$. Finally, $f(2)=(\sqrt2-1)^2>0$.
Since $4<(8/5)^3$, we have $a_3=\sqrt[3]{4}<8/5$. The function
$x^2-x-1$ is increasing for $x>1/2$, so
\[
 f(3)=2(a_3^2-a_3-1)
 <2\left(\frac{64}{25}-\frac85-1\right)=-\frac2{25}<0.
\]
Thus $f$ is strictly decreasing and changes sign exactly once between $2$
and $3$.

\subsection{Rational phase margins}\label{app:rational-margins}

\begin{proof}[Proof of Lemma~\ref{lem:phase-rational-margins}]
Recall from Theorem~\ref{thm:contact-phase} that
$f(p)=2a_p(a_p-1)-(p-1)$, where $a_p=2^{1-1/p}$, is strictly
decreasing and vanishes exactly at $p_c$. The exact power
comparisons
\[
 37^8<32\cdot24^8,\qquad
 128\cdot25^{11}<39^{11}
\]
give $2^{5/8}>37/24$ and $2^{7/11}<39/25$.
Since $2a(a-1)$ is increasing for $a>1/2$,
\[
 f(8/3)>2\frac{37}{24}\frac{13}{24}-\frac53
 =\frac1{288}>0,
\]
whereas
\[
 f(11/4)<2\frac{39}{25}\frac{14}{25}-\frac74
 =-\frac7{2500}<0.
\]
This proves the two rational bounds for $p_c$.

For $p_c\le p<\infty$, put $x_p=2^{1/p}$. The lower bound on $p_c$ and the first
power comparison give
\[
 x_p\le x_{p_c}<2^{3/8}<\frac{48}{37}.
\]
Corollary~\ref{cor:phase-explicit} therefore yields
\[
 \sigma_p=\frac{2-x_p}{4-x_p}>\frac{13}{50},
 \qquad
 \rho_p>\frac{13}{50}\frac{37}{50}
 =\frac{481}{2500}>\frac{19}{100}.
\]

Let $x_c=2^{1/p_c}$.  The upper bound already proved gives
\(x_c<48/37<13/10\).  On the other hand, \(p_c<11/4\) and $9^{11}<16\cdot7^{11}$ give
$9/7<2^{4/11}<x_c$. Hence $9/7<x_c<13/10$.
Using the central-contact formula,
\[
 \gamma_{p_c}-c_*=
 \frac{P(x_c)}{10000(4-x_c)^2},\qquad
 P(x)=-5623x^2+14984x-9968.
\]
The polynomial $P$ is concave, while
\[
 P(9/7)=\frac{97}{49},\qquad P(13/10)=\frac{833}{100}.
\]
Hence
\[
 \gamma_{p_c}-c_*>
 \frac{97}{3610000}>10^{-5}.
\]

Finally put $x=2^{1/4}$.  Taking fourth powers gives
$7/6<x<6/5$.  Since $s(x)=(2-x)/(4-x)$ is decreasing,
\[
 \frac27=s(6/5)<\sigma_4=s(x)<s(7/6)=\frac5{17}.
\]
The function $s(1-s)$ is increasing for $s<1/2$, and hence
$\rho_4>10/49$. The central formula
$\gamma_4=x(2-x)/(2(4-x)^2)$ is increasing in $x<4/3$, so
\[
 \gamma_4>\frac{35}{578},\qquad
 \frac{35}{578}-\frac3{50}=\frac8{14450}>10^{-5}.
\]
\end{proof}

\subsection{Critical scaling}\label{app:critical-scaling}

\begin{proposition}\label{prop:critical-scaling}
Let \(p_c\) and \(t_p\) be as in Theorem~\ref{thm:contact-phase}.
There are explicit constants $C_t,C_\kappa>0$ such that, as
$p\to p_c^-$,
\[
 t_p=C_t\sqrt{p_c-p}+o(\sqrt{p_c-p}),
\]
and
\[
 \kappa_p-\bigl(2^{1-1/p}-1\bigr)
 =C_\kappa(p_c-p)^2+o\bigl((p_c-p)^2\bigr).
\]
The constants are given explicitly in the proof.
Moreover, \(\kappa_p,\sigma_p,\gamma_p\) are \(C^1\) but not \(C^2\) at \(p_c\).
\end{proposition}

\begin{proof}
Put
\[
 a_c=2^{1-1/p_c},\qquad
 B_c=\frac{(p_c-1)(p_c-3)(2p_c-1)}{12a_c^3},
\]
and
\[
 h_c'=
 \frac{1}{2a_c}\left(
 \frac{2a_c(2a_c-1)\log2}{p_c^2}-1
 \right).
\]
The constants in the statement are
\[
 C_t=\sqrt{\frac{h_c'}{B_c}},\qquad
 C_\kappa=-\frac{(h_c')^2}{2B_c}.
\]
Since $2<p_c<3$, we have $B_c<0$. Also
$h_c'=f'(p_c)/(2a_c)<0$ by Appendix~\ref{app:center-sign}, so both
constants are positive.  The defining equation \eqref{eq:Up} is real-analytic near
$(p,t,u)=(p_c,0,a_c)$ and even in $t$. Since its derivative with respect
to $u$ is nonzero there, the analytic implicit function theorem gives a
jointly real-analytic branch $U_p(t)$ in $(p,t)$, which is even in $t$.
Consequently $W_p(z)=U_p(\sqrt z)$ is jointly analytic in $(p,z)$ near
$(p_c,0)$. Write $W_p(z)=a_p+b_pz+c_pz^2+O(z^3)$.
A second-order expansion of the defining equation gives
\[
 b_p=-\frac{p-1}{2a_p},\qquad
 c_p=\frac{(p-1)(p-3)(2p-1)}{24a_p^3},
\]
and therefore
\[
 W_p''(0)=\frac{(p-1)(p-3)(2p-1)}{12a_p^3}.
\]
Put $h(p)=\Phi_p(0)=f(p)/(2a_p)$. Since $f(p_c)=0$,
$h'(p_c)=h_c'<0$, while $\partial_z\Phi_{p_c}(0)=W_{p_c}''(0)=B_c<0$.
The implicit function theorem applied to $\Phi_p(z)=0$ gives
\[
 z_p=\frac{h_c'}{B_c}(p_c-p)+o(p_c-p).
\]
The coefficient is positive, and taking square roots gives the formula for
$t_p$. Finally, uniformly along this branch,
\[
 S_p(z)-S_p(0)
 =h(p)z+\frac{B_c}{2}z^2
  +o\bigl((|p-p_c|+z)z\bigr).
\]
Substitution of the expansion for $z_p$ gives the stated quadratic gain and
the constant $C_\kappa$.

The analytic stationary branch constructed above shows that the interior
and central branches of $\kappa_p$ agree to first order at $p_c$, while
their second derivatives differ by $2C_\kappa>0$. The formulas in \eqref{eq:intro-parameters} therefore imply
that $\kappa_p,\sigma_p,\gamma_p$ are $C^1$ but not $C^2$ at $p_c$.

\end{proof}

\section{Exact inequalities for the interior-contact argument}
\label{app:interior-scalar}

The calculations below complete the proofs of
Lemmas~\ref{lem:canonical-localization}
and~\ref{lem:canonical-dichotomy}.
We also give an exact certificate for the
sampling obstruction in
Example~\ref{ex:sampling-obstruction}.

\subsection{Localization inequalities}
\label{app:canonical-localization-proof}

\begin{proof}[Completion of the proof of Lemma~\ref{lem:canonical-localization}]
Write \(b=1-r=(1+q)/2\) and \(m=z+w\), and use
the mass identities~\eqref{eq:canonical-masses}. In particular, $0\le\omega\le\omega_{\max}:=(1+q)^2/16$.
The definition of $C_{\rm rnd}$ gives $C_{\rm rnd}=m/(2\sqrt2)$.
For the fixed-$I$ sample, $\lambda=q/(1+q)$ and $\mu=1/(1+q)$.
Substitution in \eqref{eq:canonical-costs} gives the main coordinates
\begin{equation}\label{eq:canonical-endpoint}
 \begin{aligned}
 X&=\frac{q^2z}{(1+q)^2}+\frac{qw}{1+q}
   =q\left(\frac{2-5q}{16}
       +\frac{2+q}{(1+q)^2}\omega-\frac{\delta_2}{1+q}\right),\\
 Y&=\frac{z}{(1+q)^2}=\frac1{16}-\frac{\omega}{(1+q)^2},
 \qquad C_{\rm fix}^2=X^2+Y^2.
 \end{aligned}
\end{equation}
Thus $X,Y\ge0$ at every feasible profile.

The main-text proof gives $m>22/125$, $q<37/125$,
$s\ge1/2-r$, and $s,\tau>0$. The endpoint estimates below give the
remaining defect bounds.

By \eqref{eq:canonical-masses},
\begin{equation}\label{eq:canonical-omega-zero}
 \omega>\omega_0(q)+\delta_2,\qquad
 \omega_0(q)=\frac{625q^2+250q-23}{2000}.
\end{equation}
For fixed $q,\omega$, replacing $\delta_2$ by zero in
\eqref{eq:canonical-endpoint} increases $X$ by
$q\delta_2/(1+q)\ge0$ and leaves $Y$ unchanged.  Since the actual $X\ge0$,
this can only increase the squared norm.  This is only an algebraic upper
envelope; the modified quadruple is not asserted to remain a feasible
graph profile.  Denote this squared envelope by $F_q(\omega)$.
It is a convex quadratic in $\omega$, so its maximum over any compact
interval is attained at an endpoint.

Suppose $q\ge1/12$.  By \eqref{eq:canonical-omega-zero}, it suffices to
check $\omega=\omega_0(q)$ and $\omega=\omega_{\max}$.
The upper endpoint gives
\[
 \sqrt{F_q(\omega_{\max})}=\frac{q(1-q)}4
 <\frac{21}{400}<c_*,
\]
since $0\le q<37/125<3/10$.  At the lower endpoint the coordinates are
\[
 X_0=\frac{q(125q^2+88q+51)}{500(1+q)^2},\qquad
 Y_0=\frac{37-125q^2}{500(1+q)^2}.
\]
Write
\[
 D_0(q)=Y_0-X_0=\frac{37-125q}{500},\qquad
 S_0(q)=X_0+Y_0
 =\frac{125q^3-37q^2+51q+37}{500(1+q)^2}.
\]
Then
\[
 S_0'(q)=\frac{P(q)}{500(1+q)^3},\qquad
 P(q)=125q^3+375q^2-125q-23.
\]
The numerator is convex on $[1/12,3/10]$, since
$P''(q)=750q+750>0$, and
\[
 P(1/12)=-\frac{53119}{1728}<0,\qquad P(3/10)=-\frac{187}{8}<0.
\]
Hence $S_0'<0$.  On the relevant interval $[1/12,37/125)$,
$S_0,D_0$ are nonnegative and decreasing.  Since
\[
 S_0(1/12)=\frac{70961}{1014000}<\frac7{100},\qquad
 D_0(1/12)=\frac{319}{6000}<\frac4{75},
\]
we obtain
\[
 F_q(\omega_0(q))=\frac{S_0(q)^2+D_0(q)^2}{2}
 <\frac{697}{180000}<c_*^2,
 \qquad c_*^2-\frac{697}{180000}=\frac{8161}{900000000}>0.
\]
Both endpoints contradict $C_{\rm fix}>c_*$, proving $q<1/12$.

Next suppose $\omega\ge\varepsilon:=1/1000$.  The same convexity
argument reduces to $\omega=\varepsilon$ and $\omega=\omega_{\max}$;
the latter was just handled.  At $\omega=\varepsilon$, the envelope
coordinates are
\[
 X_\varepsilon(q)=\frac{q(2-5q)}{16}
 +\frac1{1000}\left(1-\frac1{(1+q)^2}\right),\qquad
 Y_\varepsilon(q)=\frac1{16}-\frac1{1000(1+q)^2}.
\]
Both are nonnegative and increasing on $[0,1/12]$.  At $q=1/12$,
\[
 X_\varepsilon=\frac{16343}{1946880}<\frac{17}{2000},\qquad
 Y_\varepsilon=\frac{20837}{338000}<\frac{617}{10000}.
\]
Thus $F_q(\varepsilon)<c_*^2$, with the rational certificate
\[
 c_*^2-\left[\left(\frac{17}{2000}\right)^2
 +\left(\frac{617}{10000}\right)^2\right]
 =\frac{43}{20000000}>0.
\]
This contradicts $C_{\rm fix}>c_*$, so $\omega<1/1000$.
Finally \eqref{eq:canonical-masses} yields
\[
 \delta_2<\frac{3-2q-5q^2}{16}+\frac1{1000}-\frac{22}{125}
 =\frac1{80}-\frac q8-\frac{5q^2}{16}.
\]
Since $\delta_2\ge0$, the actual double-failure branch also has $h(q)>0$.
\end{proof}

\subsection{The canonical envelope estimate}
\label{app:canonical-dichotomy-proof}

We verify the monotonicity and endpoint estimate used in the proof of
Lemma~\ref{lem:canonical-dichotomy}.  Here
$\varepsilon=1/1000$, $\rho=\sqrt{10}/50$,
$b=(1+q)/2$, $r=(1-q)/2$, and
$R_*(q)=(b+\rho)/(b-\rho)$.  Let $A(q)$ denote the second term in
\eqref{eq:canonical-envelope}.  Then
\[
 \frac{A'(q)}{A(q)}=-\frac{q}{2br}
 -\frac{b}{2(b^2-4\varepsilon)}<0.
\]
For
\[
 G(q)=R_*(q)h(q)+\frac q4(b+\rho),
\]
the relations $b' =1/2$ and $h(q)=1/80-q/8-5q^2/16$ give
\[
 G'(q)=
 -\frac{\rho(3/10-q^2-\rho^2)}{4(b-\rho)^2}
 -\frac{q(3b+5\rho)}{8(b-\rho)}<0.
\]
Indeed $b>\rho$, and
$q^2+\rho^2<1/144+1/250<3/10$ on $0\le q<1/12$.
Thus $H=A+G+\varepsilon$ is decreasing on $0\le q<1/12$.  Finally,
\[
 H(0)=\frac1{1000}+\frac1{2\sqrt{123}}
 +\frac{5\sqrt{10}+2}{5\sqrt{10}-2}\frac1{80}.
\]
The elementary estimates
\[
 \frac1{2\sqrt{123}}<\frac7{155},\qquad
 \frac{5\sqrt{10}+2}{5\sqrt{10}-2}<\frac{129}{100}
\]
follow from $155^2<14^2\cdot123$ and $458^2<145^2\cdot10$, respectively.
Hence
\[
 H(0)<\frac1{1000}+\frac7{155}+\frac{129}{8000}
 =\frac{15447}{248000},\qquad
 c_*-H(0)>\frac{17}{1240000}>0.
\]
This proves the required bound.

\subsection{Exact certificate for the sampling obstruction}
\label{app:full-sampling-obstruction}

We verify the claim in Example~\ref{ex:sampling-obstruction} directly for
its concrete profile. Consider the main terms from a balanced two-sided
sample: a proportion $\theta$ of $J$ and a proportion $\xi$ of $I$ are
placed in the first class. Since $r=10/21$, balance gives
\[
 \frac{10}{21}\xi+\frac{11}{21}\theta=\frac12,
 \qquad
 \xi=\frac{21}{20}-\frac{11}{10}\theta.
\]
The feasibility condition $0\le\xi\le1$ is equivalent to
$1/22\le\theta\le21/22$. Substituting
$z=121/1764$ and $w=11/98$ into the expected normalized internal-edge
counts gives
\[
 a_\theta=\theta^2z+\xi\theta w
 =\frac{11\theta(189-88\theta)}{17640},
\]
\[
 b_\theta=(1-\theta)^2z+(1-\xi)(1-\theta)w
 =\frac{11(1-\theta)(101+88\theta)}{17640}.
\]
Put $u=(1-2\theta)^2$. As $\theta$ ranges over its feasible interval,
$u$ ranges over $[0,100/121]$, and direct expansion yields
\[
 a_\theta^2+b_\theta^2
 =\frac{14641}{38896200}
   \left(u-\frac{2559}{3872}\right)^2
 +\frac{156269119}{39829708800}.
\]
Since $2559/3872\in(0,100/121)$, the minimum is therefore
\[
 \min_\theta(a_\theta^2+b_\theta^2)
 =\frac{156269119}{39829708800}
 =\frac1{256}+\frac{684319}{39829708800}>\frac1{256}.
\]
Thus every main-term vector in this two-sided sampling family has
normalized $\ell_2$ cost strictly larger than $1/16$.

To extend the obstruction to nearby norms, let
\[
 m(p)=\min_{1/22\le\theta\le21/22}\|(a_\theta,b_\theta)\|_p.
\]
The parameter interval is compact, so $m(p)$ depends continuously on $p$,
while $\gamma_p\to1/16$ as $p\downarrow2$. Since the strict inequality above
gives $m(2)>1/16$, we have $m(p)>\gamma_p$ for every $p>2$ sufficiently
close to $2$. This is the persistence assertion used in
Example~\ref{ex:sampling-obstruction}.

\section*{Declaration on the use of AI}
The author used generative AI tools to assist in discussing proof strategies, checking proofs, and improving exposition.

\end{document}